\documentclass[11pt,notitlepage]{article}
\usepackage{amssymb, amsmath, amsthm}

\usepackage{amsfonts,mathtools,mathrsfs,mathtools}
\usepackage{color}
\usepackage[colorlinks,linkcolor=blue,citecolor=red]{hyperref}
\usepackage{graphicx,caption,subcaption}

\newtheorem{theorem}{Theorem}[section]
\newtheorem{lemma}[theorem]{Lemma}

\newtheorem{corollary}[theorem]{Corollary}

\numberwithin{equation}{section}

\newcommand{\R}{\mathbb{R}}  
\renewcommand{\L}{\mathcal{L}}  
\renewcommand{\H}{\mathcal{H}}  
\newcommand{\X}{\mathcal{X}}  
\newcommand{\E}{\mathcal{E}}
\newcommand{\T}{\mathbb{T}}

\begin{document}

\title{A Variational Perspective on Accelerated Methods\\ in Optimization}
\author{Andre Wibisono \\ \texttt{wibisono@berkeley.edu}
\and
Ashia C. Wilson \\ \texttt{ashia@berkeley.edu}
\and
Michael I. Jordan \\ \texttt{jordan@cs.berkeley.edu}}
\maketitle

\begin{abstract}
Accelerated gradient methods play a central role in optimization, achieving
optimal rates in many settings.  While many generalizations and extensions
of Nesterov's original acceleration method have been proposed, it is not
yet clear what is the natural scope of the acceleration concept.  
In this paper, we study accelerated methods from a continuous-time perspective.
We show that there is a Lagrangian functional that we call the \emph{Bregman
Lagrangian} which generates a large class of accelerated methods in continuous
time, including (but not limited to) accelerated gradient descent, its
non-Euclidean extension, and accelerated higher-order gradient methods.
We show that the continuous-time limit of all of these methods correspond
to traveling the same curve in spacetime at different speeds.  From this
perspective, Nesterov's technique and many of its generalizations can be
viewed as a systematic way to go from the continuous-time curves generated
by the Bregman Lagrangian to a family of discrete-time accelerated algorithms. 
\end{abstract}

\section{Introduction}
\label{Sec:Intro}

The phenomenon of \emph{acceleration} plays an important role in theory and practice of
convex optimization.  Introduced by Nesterov in 1983~\cite{Nesterov83} in the context of
gradient descent,  the acceleration idea has been extended to a wide
range of other settings, including composite optimization~\cite{Nesterov07,Tseng08,BeckTeboulle09}, 
stochastic optimization~\cite{HuKwokPan09,Lan12}, nonconvex optimization~\cite{GhadimiLan15,LiLin15},
and conic programming~\cite{LanLuMonteiro11}.  There have been generalizations
to non-Euclidean optimization~\cite{Nesterov05,Krichene15} and higher-order algorithms~\cite{Nesterov08,
Baes09}, and there have been numerous applications that further extend the reach
of the idea~\cite{Ji09,JiSJY09,JGK10,Mukherjee13}.  On the theoretical front,
acceleration often improves the convergence rate of the underlying gradient-based
procedure, and, under certain conditions, yields an optimal rate~\cite{NemirovskiiYudin}.

Despite this compelling evidence of the value of the idea of acceleration, it remains something
of a conceptual mystery.  Derivations of accelerated methods do not flow from a single
underlying principle, but tend to rely on case-specific algebra~\cite{Juditsky13}.  The basic
Nesterov technique is often explained intuitively in terms of momentum, but this intuition
does not easily carry over to non-Euclidean settings~\cite{Allen-Zhu14}.  In recent years, the
number of explanations and interpretations of acceleration has increased~\cite{Allen-Zhu14,
Arjevani15,Flammarion15, Lessard14, BubeckLeeSingh15}, but these explanations have been
focused on restrictive instances of acceleration, such as first-order algorithms,
the Euclidean setting, or cases in which the objective function is strongly convex or quadratic.
It is not yet clear what the natural scope of the acceleration concept is and indeed whether
it is a single phenomenon.

In this paper we study acceleration from a continuous-time, variational point of view.
We build on recent work by~\cite{SuBoydCandes14}, who show that the continuous-time limit
of Nesterov's accelerated gradient descent is a second-order differential equation, and we take
inspiration from continuous time analysis of mirror descent~\cite{NemirovskiiYudin}.
In our approach, rather than starting from existing discrete-time accelerated gradient methods
and deriving differential equations by taking limits, we take as our point of departure a
variational formulation in which we define a functional on continuous-time curves that we refer
to as a \emph{Bregman Lagrangian}.  Next, we calculate and discretize the Euler-Lagrange equation
corresponding to the Bregman Lagrangian. It turns out that naive discretization (the Euler method)
does not yield a stable discrete-time algorithm that retains the rate of the underlying differential
equation; rather, a more elaborate discretization involving an auxiliary sequence is necessary.
This auxiliary sequence is essentially that used by Nesterov in his constructions of accelerated mirror
descent~\cite{Nesterov05} and accelerated cubic-regularized Newton's method~\cite{Nesterov08}, 
and later generalized by Baes~\cite{Baes09}.
Thus, from our perspective, Nesterov's approach can be viewed as a methodology for
the discretization of a certain class of differential equations.  Given the complexities
associated with the discretization of differential equations, it is perhaps not surprising
that it has been difficult to perceive the generality and scope of the acceleration concept
in a discrete-time framework.

The Bregman-Lagrangian framework permits a systematic understanding of the matching rates
associated with higher-order gradient methods in discrete and continuous time.  In the
case of gradient descent, Su et al.\ show that the discrete and continuous-time dynamics
have convergence rates of $O(1/(\epsilon k))$ and $O(1/t)$, respectively, and that these
match using the identification $t = \epsilon k$; for accelerated gradient descent,
the convergence rates are $O(1/(\epsilon k^2))$ and $O(1/t^2)$ respectively, which match
using the identification $t = \sqrt{\epsilon} k$~\cite{SuBoydCandes14}.
This result has been extended to the non-Euclidean case by Krichene et al.~\cite{Krichene15}. Higher-order
gradient descent is a descent method which minimizes a regularized $(p-1)$-st order
Taylor approximation of the objective function $f$, generalizing gradient descent
($p=2$) and Nesterov and Polyak's cubic-regularized Newton's method ($p=3$)~\cite{NesterovPolyak06}.
The $p$-th order gradient algorithm with a constant step size $\epsilon$ has convergence
rate $O(1/(\epsilon k^{p-1}))$ when $\nabla^{p-1} f$ is $(1/\epsilon)$-Lipchitz
and, in continuous time, as $\epsilon \to 0$, this algorithm corresponds to the $p$-th
\emph{rescaled gradient flow}, which is a first-order differential equation with a matching
convergence rate $O(1/t^{p-1})$.  Thus, the $p$-th order gradient algorithm can be seen
as a discretization $t = \delta k$ of the rescaled gradient flow with time step
$\delta = \epsilon^{1/(p-1)}$.  Similarly, we show that the accelerated higher-order gradient
algorithm achieves an improved convergence rate $O(1/(\epsilon k^p))$ under the same
assumption (i.e., $\nabla^{p-1} f$ is $(1/\epsilon)$-Lipschitz).  In continuous time,
as $\epsilon \to 0$, this corresponds to the second-order Euler-Lagrange curve of the Bregman Lagrangian
with a matching convergence rate $O(1/t^p)$.  Thus, the $p$-th order accelerated algorithm
can be seen as a discretization $t = \delta k$ of the Euler-Lagrange equation of the
Bregman Lagrangian with time step $\delta = \epsilon^{1/p}$.

In addition to its value in relating continuous-time and discrete-time acceleration,
the study of the Bregman Lagrangian can provide further insights into the nature of
acceleration.  For instance, it is noteworthy that the Bregman Lagrangian is closed
under time dilation.  This means that if we take an Euler-Lagrange curve of a Bregman
Lagrangian and reparameterize time so we travel the curve at a different speed, then
the resulting curve is also the Euler-Lagrange curve of another Bregman Lagrangian,
with appropriately modified parameters.  Thus, the entire family of accelerated methods
correspond to a single curve in spacetime and can be obtained by speeding up (or slowing down) any single
curve.  Another insight is obtained by noting that from the discrete-time point of
view, an interpretation of acceleration starts with a base algorithm, which we can accelerate
by coupling with a suitably weighted mirror descent step.  From the continuous-time
point of view, however, it is the weighted mirror descent step that is important since
the base gradient algorithm operates on a smaller time scale. Thus, Nesterov's accelerated
gradient methods are but one possible implementation of second-order Bregman-Lagrangian
curves as a discrete-time algorithm.

The remainder of the paper is organized as follows. In Section~\ref{Sec:BregLag}, we introduce
the general family of Bregman Lagrangians and study its properties. In Section~\ref{Sec:Poly},
we demonstrate how to discretize the Euler-Lagrange equations corresponding to the polynomial
subfamily of Bregman Lagrangians to obtain discrete-time accelerated algorithms. In particular,
we introduce the family of higher-order gradient methods which can be used to complete the
discretization. In Section~\ref{Sec:FurtherProp}, we discuss additional properties of the
Bregman Lagrangian, including gauge-invariance properties, connection to classical gradient
flows, and the correspondence with a functional that we refer to as a Bregman Hamiltonian.
Finally, we end in Section~\ref{Sec:Discussion} with a brief discussion.

\subsection{Problem setting}

We consider the optimization problem
\begin{align*}
\min_{x \in \X} \; f(x),
\end{align*}
where $\X \subseteq \R^d$ is a convex set and $f \colon \X \to \R$ is a continuously differentiable convex function. To simplify the presentation in this paper we focus on the case $\X = \R^d$. We also assume $f$ has a unique minimizer, $x^* \in \X$, satisfying the optimality condition $\nabla f(x^*) = 0$. We use the inner product norm $\|x\| = \langle x,x \rangle^{1/2}$.

We consider the general non-Euclidean setting in which the space $\X$ is endowed with a distance-generating function $h \colon \X \to \R$ that is convex and essentially smooth (i.e., $h$ is continuously differentiable in $\X$, and $\|\nabla h(x)\|_* \to \infty$ as $\|x\| \to \infty$). The function $h$ can also be used to define an alternative measure of distance in $\X$ via its Bregman divergence:
\begin{align*}
D_h(y,x) = h(y) - h(x) - \langle \nabla h(x), y-x \rangle,
\end{align*}
which is nonnegative since $h$ is convex. 
When $x$ and $y$ are nearby the Bregman divergence is an approximation to the Hessian metric,
\begin{align*}
D_h(y,x) \approx \frac{1}{2} \langle y-x, \nabla^2 h(x) (y-x) \rangle := \frac{1}{2} \|y-x\|^2_{\nabla^2 h(x)}.
\end{align*}
The \emph{Euclidean setting} is obtained when $h(x) = \frac{1}{2} \|x\|^2$, in which case the Bregman divergence and Hessian metric coincide since $\nabla^2 h(x)$ is the identity matrix.

In continuous time, the Hessian metric is generally studied rather than the more general Bregman divergence; this is the case, for instance, in the case of natural gradient flow, which is the continuous-time limit of mirror descent \cite{Alvarez,Raskutti15}. By way of contrast, we shall see that our continuous-time, Lagrangian framework crucially employs the Bregman divergence.

In this paper we denote a discrete-time sequence in lower case, e.g., $x_k$ with $k \ge 0$ an integer. We denote a continuous-time curve in upper case, e.g., $X_t$ with $t \in \R$. An over-dot means derivative with respect to time, i.e., $\dot X_t = \frac{d}{dt} X_t$.

\section{The Bregman Lagrangian}
\label{Sec:BregLag}

We define the \emph{Bregman Lagrangian}
\begin{equation}\label{Eq:BregLag}
\L(X, V, t) = e^{\alpha_t + \gamma_t}\left(D_h(X+ e^{-\alpha_t} V, X) - e^{\beta_t} f(X)\right) 
\end{equation}
which is a function of position $X \in \X$, velocity $V \in \R^d$, and time $t \in \T$, where $\T \subseteq \R$ is an interval of time. The functions $\alpha, \beta, \gamma \colon \T \to \R$ are arbitrary smooth (continuously differentiable) functions of time that determine the weighting of the velocity, the potential function, and the overall damping of the Lagrangian.
We also define the following \emph{ideal scaling conditions}:
\begin{subequations}\label{Eq:IdeSca}
\begin{align}
\dot \beta_t \,&\leq\, e^{\alpha_t}   \label{Eq:IdeScaBet} \\
\dot \gamma_t \,&=\, e^{\alpha_t}  \label{Eq:IdeScaGam};
\end{align}
\end{subequations}
these conditions will be justified in the following section.

\subsection{Convergence rates of the Euler-Lagrange equation}
\label{Sec:BregLagRate}

In this section we show that---under the ideal scaling assumption~\eqref{Eq:IdeSca}---the Bregman Lagrangian~\eqref{Eq:BregLag} defines a variational problem the solutions to which minimize the objective function $f$ at an exponential rate.

Given a general Lagrangian $\L(X, V, t)$, we define a functional on curves $\{X_t : t \in \T\}$ via integration of the Lagrangian: $J(X) = \int_{\T} \L(X_t, \dot X_t, t) dt$.  From the calculus of variations, a necessary condition for a curve to minimize this functional is that it solve the Euler-Lagrange equation:
\begin{align}\label{Eq:EL}
\frac{d}{dt}\left\{\frac{\partial \L}{\partial V} (X_t, \dot X_t, t)\right\} = \frac{\partial \L}{\partial X} (X_t, \dot X_t, t).
\end{align}
Specifically, for the Bregman Lagrangian~\eqref{Eq:BregLag}, the partial derivatives are
\begin{subequations}\label{Eq:BregLagPartial}
\begin{align}
\frac{\partial \L}{\partial X}(X, V, t)
    \,&=\, e^{\gamma_t + \alpha_t} \left( \nabla h (X + e^{-\alpha_t} V) - \nabla h(X) - e^{-\alpha_t} \nabla^2 h(X) \, V - e^{\beta_t} \nabla f(X) \right) \\
\frac{\partial \L}{\partial V}(X, V, t)
    \,&=\, e^{\gamma_t} \left( \nabla h (X + e^{-\alpha_t} V) - \nabla h(X)  \right).
\end{align}
\end{subequations}
Thus, for general functions $\alpha_t, \beta_t, \gamma_t$, the Euler-Lagrange equation~\eqref{Eq:EL} for the Bregman Lagrangian~\eqref{Eq:BregLag} is a second-order differential equation given by
\begin{equation}
\begin{split} \label{Eq:ELBreg} 
 \ddot X_t &+ (e^{\alpha_t} - \dot \alpha_t) \dot X_t + e^{2\alpha_t + \beta_t}\Big[\nabla^2 h(X_t + e^{-\alpha_t}\dot X_t)\Big]^{-1}\nabla f(X_t) \\
& + e^{\alpha_t}(\dot \gamma_t - e^{\alpha_t})\Big[\nabla^2 h(X_t + e^{-\alpha_t}\dot X_t)\Big]^{-1}(\nabla h(X_t + e^{-\alpha_t} \dot X_t) - \nabla h(X_t)) = 0.
\end{split}
\end{equation}

We now impose the ideal scaling condition~\eqref{Eq:IdeScaGam}.  In this case the last term in~\eqref{Eq:ELBreg} vanishes, so the Euler-Lagrange equation simplifies to 
\begin{align} \label{Eq:ELBreg2} 
\ddot X_t + (e^{\alpha_t} - \dot \alpha_t) \dot X_t + e^{2\alpha_t + \beta_t}\Big[\nabla^2 h(X_t + e^{-\alpha_t}\dot X_t)\Big]^{-1}\nabla f(X_t) = 0.
\end{align}
In~\eqref{Eq:ELBreg2}, we have assumed the Hessian matrix $\nabla^2 h(X_t + e^{-\alpha_t}\dot X_t)$ is invertible. But we can also write the equation~\eqref{Eq:ELBreg2} in the following way, which only requires that $\nabla h$ be differentiable,
\begin{align} \label{Eq:ELBreg3} 
\frac{d}{dt} \nabla h(X_t + e^{-\alpha_t} \dot X_t) = -e^{\alpha_t+\beta_t} \nabla f(X_t).
\end{align}

To establish a convergence rate associated with solutions to the Euler-Lagrange equation---under the ideal scaling conditions---we take a Lyapunov function approach.  Defining the following energy functional:
\begin{align}\label{Eq:E}
\E_t \,=\, D_h\left(x^*, \, X_t + e^{-\alpha_t} \dot X_t\right) + e^{\beta_t} (f(X_t)-f(x^*)),
\end{align}
we immediately obtain a convergence rate, as shown in the following theorem.

\begin{theorem}\label{Thm:Rate}
If the ideal scaling~\eqref{Eq:IdeSca} holds, then solutions to the Euler-Lagrange equation~\eqref{Eq:ELBreg3} satisfy
\begin{align*}
f(X_t) - f(x^*) \le O(e^{-\beta_t}).
\end{align*}
\end{theorem}
\begin{proof}
The time derivative of the energy functional is
\begin{align*}
\dot \E_t = -\left\langle \frac{d}{dt} \nabla h(X_t + e^{-\alpha_t} \dot X_t), \, x^* - X_t - e^{-\alpha_t} \dot X_t \right \rangle + \dot \beta_t e^{\beta_t} (f(X_t) - f(x^*)) + e^{\beta_t} \langle \nabla f(X_t), \dot X_t \rangle.
\end{align*}
If $X_t$ satisfies the Euler-Lagrange equation~\eqref{Eq:ELBreg3}, then the time derivative simplifies to
\begin{align*}
\dot \E_t 
&= -e^{\alpha_t+\beta_t} D_f(x^*, X_t) + (\dot \beta_t - e^{\alpha_t}) e^{\beta_t} (f(X_t) - f(x^*))
\end{align*}
where $D_f(x^*,X_t) = f(x^*) - f(X_t) - \langle \nabla f(X_t), x^* - X_t \rangle$ is the Bregman divergence of $f$.
Note that $D_f(x^*,X_t) \ge 0$ since $f$ is convex, so the first term in $\dot \E_t$ is nonpositive. Furthermore, if the ideal scaling condition~\eqref{Eq:IdeScaBet} holds, then the second term is also nonpositive, so $\dot \E_t \le 0$. Since $D_h(x^*, X_t + e^{-\alpha_t} \dot X_t) \ge 0$, this implies that for any $t \ge t_0 \in \T$, $e^{\beta_t}(f(X_t) - f(x^*)) \,\le\, \E_t \,\le\, \E_{t_0}$. Thus, $f(X_t) - f(x^*) \le \E_{t_0} e^{-\beta_t} = O(e^{-\beta_t})$, as desired.
\end{proof}

For a given $\alpha_t$, which determines $\gamma_t$ by~\eqref{Eq:IdeScaBet}, the optimal convergence rate is achieved by setting $\dot \beta_t = e^{\alpha_t}$, resulting in convergence rate $O(e^{-\beta_t}) = O(\exp(-\int_{t_0}^t e^{\alpha_s}\, ds))$.
In Section~\ref{Sec:Poly} we study a subfamily of Bregman Lagrangians that have a polynomial convergence rate, and we show how we can discretize the resulting Euler-Lagrange equations to obtain discrete-time methods that have a matching, accelerated convergence rate.
In Section~\ref{Sec:Exp} we study another subfamily of Bregman Lagrangians that have an exponential convergence rate, and discuss its connection to a generalization of Nesterov's restart scheme. 
In the Euclidean setting, our derivations simplify. We present these derivations in Appendix~\ref{App:Euc}, and comment on the insight that they provide into the question posed by Su et al.~\cite{SuBoydCandes14} on the significance of the value $3$ in the damping coefficient for Nesterov's accelerated gradient descent.

\subsection{Time dilation}
\label{Sec:BregLagTime}

A notable property of the Bregman Lagrangian family is that it is closed under time dilation. 
This means if we take the Euler-Lagrange equation~\eqref{Eq:ELBreg} of the Bregman Lagrangian~\eqref{Eq:BregLag} and reparameterize time to travel the curve at a different speed, the resulting curve is also the Euler-Lagrange equation of a Bregman Lagrangian with a suitably modified set of parameters.

Concretely, let $\tau \colon \T \to \T'$ be a smooth (twice-continuously differentiable) increasing function, where $\T' = \tau(\T) \subseteq \R$ is the image of $\T$. Given a curve $X \colon \T' \to \X$, we consider the reparameterized curve $Y \colon \T \to \X$ defined by
\begin{align}\label{Eq:Reparam}
Y_t = X_{\tau(t)}.
\end{align}
That is, the new curve $Y$ is obtained by traversing the original curve $X$ at a new speed of time determined by $\tau$. If $\tau(t) > t$, then we say that $Y$ is the {\em sped-up version} of $X$, because the curve $Y$ at time $t$ has the same value as the original curve $X$ at the future time $\tau(t)$.

For clarity, we let $\L_{\alpha,\beta,\gamma}$ denote the Bregman Lagrangian~\eqref{Eq:BregLag} parameterized by $\alpha, \beta, \gamma$. Then we have the following result whose proof is provided in Appendix~\ref{App:TimeProof}.

\begin{theorem}\label{Thm:Time}
If $X_t$ satisfies the Euler-Lagrange equation~\eqref{Eq:ELBreg} for the Bregman Lagrangian $\L_{\alpha,\beta,\gamma}$, then the reparameterized curve $Y_t = X_{\tau(t)}$ satisfies the Euler-Lagrange equation for the Bregman Lagrangian $\L_{\tilde \alpha, \tilde \beta, \tilde \gamma}$, with modified parameters 
\begin{subequations}\label{Eq:ABC}
\begin{align}
\tilde \alpha_t &= \alpha_{\tau(t)} + \log \dot \tau(t) \\
\tilde \beta_t &= \beta_{\tau(t)} \\
\tilde \gamma_t &= \gamma_{\tau(t)}.
\end{align}
\end{subequations}
Furthermore, $\alpha,\beta,\gamma$ satisfy the ideal scaling~\eqref{Eq:IdeSca} if and only if $\tilde \alpha, \tilde \beta, \tilde \gamma$ do.
\end{theorem}

We note that in general, when we reparameterize time by a time-dilation function $\tau(t)$, the Lagrangian functional transforms to $\tilde \L(X, V, t) = \dot \tau(t) \, \L\left(X, \frac{1}{\dot \tau(t)} V, \tau(t) \right).$
Thus, another way of stating the result in Theorem~\ref{Thm:Time} is to claim that
\begin{align}\label{Eq:LagTime}
\L_{\tilde \alpha, \tilde \beta, \tilde \gamma}(X, V, t) = \dot \tau(t) \, \L_{\alpha,\beta,\gamma}\left(X, \frac{1}{\dot \tau(t)} V, \tau(t) \right),
\end{align}
which we can easily verify by directly substituting the definition of the Lagrangian~\eqref{Eq:BregLag} and the modified parameters $\tilde \alpha, \tilde \beta, \tilde \gamma$~\eqref{Eq:ABC}.

In Section~\ref{Sec:Poly}, we show that the Bregman Lagrangian generates the family of higher-order accelerated methods in discrete time. Thus, the time-dilation property means that the entire family of curves for accelerated methods in continuous time corresponds to a single curve in spacetime, which is traveled at different speeds. This suggests that the underlying solution curve  has a more fundamental structure that is worth exploring further.

\section{Polynomial convergence rates and accelerated methods}
\label{Sec:Poly}

In this section, we study a subfamily of Bregman Lagrangians~\eqref{Eq:BregLag} with the following choice of parameters, indexed by a parameter $p > 0$,
\begin{subequations}\label{Eq:ABCPoly}
\begin{align}
\alpha_t &=  \log p - \log t\\
\beta_t &= p\log t + \log C\\
\gamma_t &= p \log t,
\end{align}
\end{subequations}
where $C > 0$ is a constant. The parameters $\alpha,\beta,\gamma$ satisfy the ideal scaling condition~\eqref{Eq:IdeSca} (with an equality on the first condition~\eqref{Eq:IdeScaBet}). The Euler-Lagrange equation~\eqref{Eq:ELBreg2} is given by
\begin{align}\label{Eq:ELPoly}
\ddot X_t + \frac{p+1}{t} \dot X_t + Cp^2 t^{p-2} \left[\nabla^2 h\left(X_t + \frac{t}{p} \dot X_t\right)\right]^{-1} \nabla f(X_t) = 0
\end{align}
and, by Theorem~\ref{Thm:Rate}, it has an $O(1/t^p)$ rate of convergence.
As direct result of the time-dilation property (Theorem~\ref{Thm:Time}), the entire family of curves~\eqref{Eq:ELPoly} can be obtained by speeding up the curve in the case $p=2$ by the time-dilation function $\tau(t) = t^{p/2}$.
In Appendix~\ref{App:Existence} we discuss the issue of the existence and uniqueness of the solution to the differential equation~\eqref{Eq:ELPoly}.

The case $p=2$ of the equation~\eqref{Eq:ELPoly} is the continuous-time limit of Nesterov's accelerated mirror descent~\cite{Nesterov05}, and the case $p=3$ is the continuous-time limit of Nesterov's accelerated cubic-regularized Newton's method~\cite{Nesterov08}.
The case $p=2$ has also been derived independently in a recent work of Krichene et al.~\cite{Krichene15}; in the Euclidean case, when the Hessian $\nabla^2 h$ is the identity matrix, we recover the differential equation of Su et al.~\cite{SuBoydCandes14}.

\subsection{Naive discretization}

We now turn to the challenge of discretizing the differential equation in~\eqref{Eq:ELPoly}, with the goal of obtaining a
discrete-time algorithm whose convergence rate matches that of the underlying differential equation.  As we show in this
section, a naive Euler method is not able to match the underlying rate.  To match the rate a more sophisticated
approach is needed, and it is at this juncture that Nesterov's three-sequence idea makes its appearance.

We first write the second-order equation~\eqref{Eq:ELPoly} as the following system of first-order equations:
\begin{subequations}\label{Eq:ELPoly1}
 \begin{align}
 Z_t &= X_t + \frac{t}{p}\dot X_t    \label{Eq:ELPoly1a} \\
\frac{d}{dt} \nabla h(Z_t) &= -Cpt^{p-1}\nabla f(X_t)   \label{Eq:ELPoly1b}.
 \end{align}
\end{subequations}
Now we discretize $X_t$ and $Z_t$ into sequences $x_k$ and $z_k$ with time step $\delta >0$. That is, we make the identification $t = \delta k$ and set $x_k = X_t$,  $x_{k+1} = X_{t + \delta} \approx X_t + \delta \dot X_t$ and $z_k = Z_t$, $z_{k+1} = Z_{t + \delta} \approx Z_t + \delta \dot Z_t$. Applying the forward-Euler method to~\eqref{Eq:ELPoly1a} gives the equation $z_k = x_k + \frac{\delta k}{p} \frac{1}{\delta}(x_{k+1}-x_k)$, or equivalently,
\begin{align}\label{Eq:Alg0x}
x_{k+1} = \frac{p}{k} z_k + \frac{k-p}{k} x_k.
\end{align}
Similarly, applying the backward-Euler method to equation~\eqref{Eq:ELPoly1b} gives $\frac{1}{\delta}(\nabla h(z_{k}) - \nabla h(z_{k-1})) = -Cp(\delta k)^{p-1} \nabla f(x_k)$, which we can write as the optimality condition of the following weighted mirror descent step:
\begin{align}\label{Eq:Alg0z}
z_k = \arg\min_z \,\left\{ Cp k^{p-1} \langle \nabla f(x_k), z \rangle + \frac{1}{\epsilon} D_h(z,z_{k-1}) \right\},
\end{align}
with step size $\epsilon = \delta^p$. In principle, the two updates~\eqref{Eq:Alg0x},~\eqref{Eq:Alg0z} define an algorithm that implements the dynamics \eqref{Eq:ELPoly1} in discrete time. However, we cannot establish a convergence rate for the algorithm~\eqref{Eq:Alg0x},~\eqref{Eq:Alg0z}; indeed, empirically, we find that the algorithm is unstable. Even for the simple case in which $f$ is a quadratic function in two dimensions, the iterates of the algorithm initially approach and oscillate near the minimizer, but eventually the oscillation increases and the iterates shoot off to infinity.

\subsection{A rate-matching discretization}
\label{Sec:ModDisc}

We now discuss how to modify the naive discretization scheme~\eqref{Eq:Alg0x},~\eqref{Eq:Alg0z} into an algorithm whose rate matches that of the underlying differential equation. Our approach is inspired by Nesterov's constructions of accelerated mirror descent~\cite{Nesterov05} and accelerated cubic-regularized Newton's method~\cite{Nesterov08}, which maintain three sequences in the algorithms and use the estimate sequence technique to prove convergence. Indeed, from our point of view, Nesterov's methodology can be viewed as a rate-matching discretization methodology.

Specifically, we consider the following scheme, in which we introduce a third sequence $y_k$ to replace $x_k$ in the updates,
\begin{subequations}\label{Eq:Alg}
\begin{align}
x_{k+1} &= \frac{p}{k+p} z_k + \frac{k}{k+p} y_k   \label{Eq:Algx} \\
z_k &= \arg\min_z \left\{Cp k^{(p-1)} \langle \nabla f(y_k), z \rangle + \frac{1}{\epsilon} D_h(z,z_{k-1}) \right\},
\label{Eq:Algz}
\end{align}
\end{subequations}
where $k^{(p-1)} := k(k+1) \cdots (k+p-2)$ is the rising factorial. 
A sufficient condition for the algorithm~\eqref{Eq:Alg} to have an $O(1/(\epsilon k^p))$ convergence rate is that the new sequence $y_k$ satisfy the inequality
\begin{align}\label{Eq:yIneq}
\langle \nabla f(y_k), x_k - y_k \rangle \ge M \epsilon^{\frac{1}{p-1}} \|\nabla f(y_k)\|_*^{\frac{p}{p-1}},
\end{align}
for some constant $M > 0$. 
Note that in going from~\eqref{Eq:Alg0x} to~\eqref{Eq:Algx} we have replaced the weight $\frac{p}{k}$ by $\frac{p}{k+p}$; this is only for convenience in the proof given below, and does not change the asymptotics since $\frac{p}{k} = \Theta(\frac{p}{k+p})$ as $k \to \infty$. Similarly, we replace $k^{p-1}$ in~\eqref{Eq:Alg0z} by the rising factorial $k^{(p-1)}$ in~\eqref{Eq:Algz} to make the algebra easier, but we still have $k^{(p-1)} = \Theta(k^{p-1})$.

The following result also requires a uniform convexity assumption on the distance-generating function $h$. Recall that $h$ is {\em $\sigma$-uniformly convex of order $p \ge 2$} if its Bregman divergence is lower bounded by the $p$-th power of the norm,
\begin{align}\label{Eq:UnifConv}
D_h(y,x) \ge \frac{\sigma}{p} \|y-x\|^p.
\end{align}
The case $p=2$ is the usual definition of strong convexity. An example of a uniformly convex function is the $p$-th power of the norm, $h(x) = \frac{1}{p} \|x-w\|^p$ for any $w \in \X$, which is $\sigma$-uniformly convex of order $p$ with $\sigma = 2^{-p+2}$~\cite[Lemma~4]{Nesterov08}.

\begin{theorem}\label{Thm:AlgRate}
Assume $h$ is $1$-uniformly convex of order $p \ge 2$, and the sequence $y_k$ satisfies the inequality~\eqref{Eq:yIneq} for all $k \ge 0$. Then the algorithm~\eqref{Eq:Alg} with the constant $C \le M^{p-1}/p^p$ and initial condition $z_0 = x_0 \in \X$ has the convergence rate
\begin{align}\label{Eq:AlgRate}
f(y_k) - f(x^*) \le \frac{D_h(x^*,x_0)}{C \epsilon k^{(p)}} = O\left(\frac{1}{\epsilon k^p}\right).
\end{align}
\end{theorem}

The proof of Theorem~\ref{Thm:AlgRate} uses a generalization of Nesterov's estimate sequence technique, and can be found in Appendix~\ref{Thm:AlgRateProof}. We note that with the scaling $\epsilon = \delta^p$ as in the previous section, the convergence rate $O(1/(\epsilon k^p))$ matches the $O(1/t^p)$ rate in continuous time for the differential equation~\eqref{Eq:ELPoly}. We also note that the result in Theorem~\ref{Thm:AlgRate} does not require any assumptions on $f$ beyond the ability to construct a sequence $y_k$ satisfying~\eqref{Eq:yIneq}.
In the next section, we will see that we can satisfy~\eqref{Eq:yIneq} using the higher-order gradient method, which requires a higher-order smoothness assumption on $f$; the resulting algorithm is then the accelerated higher-order gradient method.

\subsection{Higher-order gradient method}

We study the higher-order gradient update, which minimizes a regularized higher-order Taylor approximation of the objective function $f$. 

Recall that for an integer $p \ge 2$, the $(p-1)$-st order Taylor approximation of $f$ centered at $x \in \X$ is the $(p-1)$-st degree polynomial
\begin{align*}
f_{p-1}(y; x) \,=\, \sum_{i=0}^{p-1} \frac{1}{i!} \nabla^i f(x) (y-x)^i
\,=\, f(x) + \langle \nabla f(x), y-x \rangle + \cdots + \frac{1}{(p-1)!} \nabla^{p-1} f(x) (y-x)^{p-1}.
\end{align*}
We say that $f$ is {\em $L$-smooth of order $p-1$} if $f$ is $p$-times continuously differentiable and 
$\nabla^{p-1} f$ is $L$-Lipschitz, which means for all $x,y \in \X$,
\begin{align}\label{Eq:Smooth}
\|\nabla^{p-1} f(y) - \nabla^{p-1} f(x)\|_* \le L \|y-x\|.
\end{align}

For a constant $N > 0$ and step size $\epsilon > 0$, we define the update operator $G_{p,\epsilon,N} \colon \X \to \X$ by
\begin{align}\label{Eq:Update}
G_{p,\epsilon,N}(x) = \arg\min_{y} \left\{ f_{p-1}(y;x) + \frac{N}{\epsilon p} \|y-x\|^p \right\}.
\end{align}
When $f$ is smooth of order $p-1$, the operator $G_{p,\epsilon,N}$ has the following property, which generalizes \cite[Lemma~6]{Nesterov08}.
We provide the proof in Appendix~\ref{App:UpdateProof}.
\begin{lemma}\label{Lem:Update}
Let $x \in \X$, $y = G_{p,\epsilon,N}(x)$, and $N > 1$. If $f$ is $L = \frac{(p-1)!}{\epsilon}$-smooth of order $p-1$, then
\begin{align}\label{Eq:UpdateIneq}
\langle \nabla f(y), x-y \rangle \,\ge\, \frac{(N^2-1)^{\frac{p-2}{2p-2}}}{2N} \, \epsilon^{\frac{1}{p-1}} \|\nabla f(y)\|_*^{\frac{p}{p-1}}.
\end{align}
Furthermore,
\begin{align}\label{Eq:UpdateIneq2}
\frac{(N^2-1)^{\frac{p-2}{2p-2}}}{2N} \, \epsilon^{\frac{1}{p-1}} \|\nabla f(y)\|_*^{\frac{1}{p-1}} \,\le\, \|x-y\| \,\le\, \frac{1}{(N-1)^{\frac{1}{p-1}}} \, \epsilon^{\frac{1}{p-1}} \|\nabla f(y)\|_*^{\frac{1}{p-1}}.
\end{align}
\end{lemma}
The inequality~\eqref{Eq:UpdateIneq} means that we can use the update operator $G_{p,\epsilon,N}$ to produce a sequence $y_k$ satisfying the requirement~\eqref{Eq:yIneq} under a higher-order order smoothness condition on $f$. We state the resulting algorithm in the next section.

\paragraph{Higher-order gradient method.}
In this section, we study the following higher-order gradient algorithm defined by the update operator $G_{p,\epsilon,N}$:
\begin{align}\label{Eq:AlgGrad}
x_{k+1} = G_{p,\epsilon,N}(x_k).
\end{align}
The case $p=2$ is the usual gradient descent algorithm, and the case $p=3$ is Nesterov and Polyak's cubic-regularized Newton's method~\cite{NesterovPolyak06}.

If $f$ is smooth of order $p-1$, then the algorithm~\eqref{Eq:AlgGrad} is a descent method. Furthermore, we can prove the following rate of convergence, which generalizes the results for gradient descent and the cubic-regularized Newton's method. We provide the proof in Appendix~\ref{App:AlgGradRateProof}.
\begin{theorem}\label{Thm:AlgGradRate}
If $f$ is $\frac{(p-1)!}{\epsilon}$-smooth of order $p-1$, then the algorithm~\eqref{Eq:AlgGrad} with constant $N > 0$ and initial condition $x_0 \in \X$ has the convergence rate
\begin{align}\label{Eq:AlgGradRate}
f(x_k) - f(x^*) \le \frac{p^{p-1} (N+1) R^p}{\epsilon k^{p-1}} = O\left(\frac{1}{\epsilon k^{p-1}}\right),
\end{align}
where $R = \sup_{x \colon f(x) \le f(x_0)} \|x-x^*\|$ is the radius of the level set of $f$ from the initial point $x_0$.
\end{theorem}

\paragraph{Rescaled gradient flow.}
We can take the continuous-time limit of the higher-order gradient algorithm as the step size $\epsilon \to 0$. The resulting curve is a first-order differential equation that is a rescaled version of gradient flow. We show that it minimizes $f$ with a matching convergence rate. 
In the following, we take $N=1$ in~\eqref{Eq:AlgGrad} for simplicity (the general $N$ simply scales the vector field by a constant). We provide the proof of Theorem~\ref{Thm:RescGradFlow} in Appendix~\ref{App:RescGradFlowProof}.

\begin{theorem}\label{Thm:RescGradFlow}
The continuous-time limit of the algorithm~\eqref{Eq:AlgGrad} is the rescaled gradient flow
\begin{align}\label{Eq:RescGradFlow}
\dot X_t = - \frac{\nabla f(X_t)}{\|\nabla f(X_t)\|_*^{\frac{p-2}{p-1}}},
\end{align}
where we define the right-hand side to be the zero if $\nabla f(X_t) = 0$. Furthermore, the rescaled gradient flow has convergence rate
\begin{align}\label{Eq:RescGradFlowRate}
f(X_t) - f(x^*) \le \frac{(p-1)^{p-1} R^p}{t^{p-1}} = O\left(\frac{1}{t^{p-1}}\right),
\end{align}
where $R = \sup_{x \colon f(x) \le f(X_0)} \|x-x^*\|$ is the radius of the level set of $f$ from the initial point $X_0$.
\end{theorem}

Equivalently, we can interpret the higher-order gradient algorithm~\eqref{Eq:AlgGrad} as a discretization of the rescaled gradient flow~\eqref{Eq:RescGradFlow} with time step $\delta = \epsilon^{\frac{1}{p-1}}$, so $t = \delta k = \epsilon^{\frac{1}{p-1}} k$. 
With this identification, the convergence rates in discrete time, $O(1/(\epsilon k^{p-1}))$, and in continuous time, $O(1/t^{p-1})$, match. The convergence rate for the continuous-time dynamics does not require any assumption beyond the convexity and differentiability of $f$ (as in the case of the Lagrangian flow~\eqref{Eq:ELBreg2}), whereas the convergence rate for the discrete-time algorithm requires the higher-order smoothness assumption on $f$.
We note that the limiting case $p \to \infty$ of~\eqref{Eq:RescGradFlow} is the {\em normalized gradient flow}, which has been shown to converge to the minimizer of $f$ in finite time~\cite{Cortes06}.
We also note that unlike the Lagrangian flow, the family of rescaled gradient flows is {\em not} closed under time dilation. 

\subsection{Accelerated higher-order gradient method}

By the result of Lemma~\ref{Lem:Update}, we see that we can use the higher-order gradient update $G_{p,\epsilon,N}$ to produce a sequence $y_k$ satisfying the inequality~\eqref{Eq:yIneq}, to complete the algorithm~\eqref{Eq:AlgGrad} that implements the polynomial family of the Bregman-Lagrangian flow~\eqref{Eq:ELPoly}. Explicitly, the resulting algorithm is as follows,
\begin{subequations}\label{Eq:AlgComp}
\begin{align}
x_{k+1} &= \frac{p}{k+p} z_k + \frac{k}{k+p} y_k   \label{Eq:AlgCompx} \\
y_k &= \arg\min_y \left\{ f_{p-1}(y;x_k) + \frac{N}{\epsilon p} \|y-x_k\|^p \right\}  \label{Eq:AlgCompy} \\
z_k &= \arg\min_z \left\{Cp k^{(p-1)} \langle \nabla f(y_k), z \rangle + \frac{1}{\epsilon} D_h(z,z_{k-1}) \right\}.  \label{Eq:AlgCompz}
\end{align}
\end{subequations}
By Theorem~\ref{Thm:AlgRate} and Lemma~\ref{Lem:Update}, we have the following guarantee for this algorithm.
\begin{corollary}\label{Cor:AlgCompRate}
Assume $f$ is $\frac{(p-1)!}{\epsilon}$-smooth of order $p-1$, and $h$ is $1$-uniformly convex of order $p$. Then the algorithm~\eqref{Eq:AlgComp} with constants $N > 1$ and $C \le (N^2-1)^{\frac{p-2}{2}}/((2N)^{p-1} p^p)$ and initial conditions $z_0 = x_0 \in \X$ has an $O(1/(\epsilon k^p))$ convergence rate.
\end{corollary}

The resulting algorithm~\eqref{Eq:AlgComp} and its convergence rate recovers the results of Baes~\cite{Baes09}, who studied a generalization of Nesterov's estimate sequence technique to higher-order algorithms. 
We note that the convergence rate $O(1/(\epsilon k^p))$ of algorithm~\eqref{Eq:AlgComp} is better than the $O(1/(\epsilon k^{p-1}))$ rate of the higher-order gradient algorithm~\eqref{Eq:AlgGrad}, under the same assumption of the $(p-1)$-st order smoothness of $f$.
This gives the interpretation of the algorithm~\eqref{Eq:AlgComp} as ``accelerating'' the higher-order gradient method. Indeed, in this view the ``base algorithm'' that we start with is the higher-order gradient algorithm in the $y$-sequence~\eqref{Eq:AlgCompy}, and the acceleration is obtained by coupling it with a suitably weighted mirror descent step in~\eqref{Eq:AlgCompx} and \eqref{Eq:AlgCompz}. 

However, from the continuous-time point of view, where our starting point is the polynomial Lagrangian flow~\eqref{Eq:ELPoly}, we see that the algorithm~\eqref{Eq:AlgComp} is only one possible implementation of the flow as a discrete-time algorithm. As we saw in Section~\ref{Sec:ModDisc}, it is only the $x$- and $z$-sequences~\eqref{Eq:AlgCompx} and~\eqref{Eq:AlgCompz} that play a role in the correspondence between the continuous-time dynamics and its discrete-time implementation, and the requirement~\eqref{Eq:yIneq} in the $y$-update is only needed to complete the convergence proof. 
Indeed, the higher-order gradient update~\eqref{Eq:AlgCompy} does not change the continuous-time limit, since from~\eqref{Eq:UpdateIneq2} in Lemma~\ref{Lem:Update} we have that $\|x_k-y_k\| = \Theta(\epsilon^{\frac{1}{p-1}})$, which is smaller than the $\delta = \epsilon^{\frac{1}{p}}$ time step in the discretization of~\eqref{Eq:ELPoly}. Therefore, the $x$ and $y$ sequences in~\eqref{Eq:AlgComp} coincide in continuous time as $\epsilon \to 0$.
Thus, from this point of view, Nesterov's accelerated methods (for the cases $p=2$ and $p=3$) are one of possibly many discretizations of the polynomial Lagrangian flow~\eqref{Eq:ELPoly}. For instance, in the case $p=2$, Krichene et al.~\cite[Section 4.1]{Krichene15} show that we can use a general regularizer in the gradient step~\eqref{Eq:AlgCompy} under some additional smoothness assumptions.
If there are other implementations, it would be interesting to see if the higher-gradient methods have some distinguishing property, such as computational efficiency.

\section{Further explorations of the Bregman Lagrangian}
\label{Sec:FurtherProp}
In addition to providing a unifying framework for the generation of accelerated gradient-based algorithms, the Bregman Lagrangian  has mathematical structure that can be investigated directly.  In this section we briefly discuss some of the additional perspective that can be obtained from the Bregman Lagrangian.  See Appendices~\ref{App:Exp}--\ref{App:NatMot} for technical details of the results discussed here.

\paragraph{Hessian vs.\ Bregman Lagrangian.}
It is important to note the presence of the Bregman divergence in the Bregman Lagrangian~\eqref{Eq:BregLag}. In the non-Euclidean setting, intuition might suggest using the Hessian metric $\nabla^2 h$ to measure a ``kinetic energy,'' and thereby obtain a \emph{Hessian Lagrangian}. This approach turns out to be unsatisfying, however, because the resulting differential equation does not yield a convergence rate and the Euler-Lagrange equation involves the third-order derivative $\nabla^3 h$, posing serious difficulties for discretization. As we have seen, the Bregman Lagrangian, on the other hand, readily provides a rate of convergence via a Lyapunov function; moreover, the resulting discrete-time algorithm in~\eqref{Eq:AlgComp} involves only the gradient $\nabla h$ via the weighted mirror descent update.
  
\paragraph{Gradient vs.\ Lagrangian flows.}
In the Euclidean case, it is known classically that we can view gradient flow as the strong-friction limit of a damped Lagrangian flow~\cite[p.~646]{Villani}. We show that the same interpretation holds for natural gradient flow and rescaled gradient flow.
In particular, we show in Appendix~\ref{App:GradLag} that we can recover natural gradient flow as the strong-friction limit of a Bregman Lagrangian flow with an appropriate choice of parameters. Similarly, we can recover the rescaled gradient flow~\eqref{Eq:RescGradFlow} as the strong-friction limit of a Lagrangian flow that uses the $p$-th power of the norm as the kinetic energy.
Therefore, the general family of second-order Lagrangian flows is more general, and includes first-order gradient flows in its closure. From this point of view, a particle with gradient-flow dynamics is operating in the regime of high friction. The particle simply rolls downhill and stops at the equilibrium point as soon as the force $-\nabla f$ vanishes; there is no oscillation since it is damped by the infinitely strong friction. Thus, the effect of moving from a first-order gradient flow to a second-order Lagrangian flow is to reduce the friction from infinity to a finite amount; this permits oscillation~\cite{ODonoghue15,SuBoydCandes14,Krichene15}, but also allows faster convergence.
   
\paragraph{Bregman Hamiltonian.}
One way to understand a Lagrangian is to study its Hamiltonian, which is the Legendre conjugate (dual function) of the Lagrangian. 
Typically, when the Lagrangian takes the form of the difference between kinetic and potential energy, the Hamiltonian is the sum of the kinetic and potential energy. The Hamiltonian is often easier to study than the Lagrangian, since its second-order Euler-Lagrangian equation is transformed into a pair of first-order equations. In our case, the Hamiltonian corresponding to the Bregman Lagrangian~\eqref{Eq:BregLag} is the following {\em Bregman Hamiltonian},
\begin{align*}
   \H(X, P, t) = e^{\alpha_t + \gamma_t} \left( D_{h^*}\left(\nabla h(X) + e^{-\gamma_t} P, \, \nabla h(X) \right) + e^{\beta_t} f(X) \right)
\end{align*}
which indeed has the form of the sum of the kinetic and potential energy. Here the kinetic energy is measured using the Bregman divergence of $h^*$, which is the convex dual function of $h$.  See Appendix~\ref{App:BregHam} for further discussion.
   
\paragraph{Gauge invariance.}
The Euler-Lagrange equation of a Lagrangian is gauge-invariant, which means it does not change when we add a total time derivative to the Lagrangian. For the Bregman Lagrangian with the ideal scaling condition~\eqref{Eq:IdeScaGam}, this property implies that we can replace the Bregman divergence $D_h(X + e^{-\alpha_t} V, X)$ in~\eqref{Eq:BregLag} by its first term $h(X + e^{-\alpha_t} V)$. This might suggest a different interpretation of the role of $h$ in the Lagrangian.

\paragraph{Natural motion.}
The \emph{natural motion} of the Bregman Lagrangian (i.e., the motion when there is no force, $-\nabla f \equiv 0$) is given by $X_t = ae^{-\gamma_t} + b$, for some constants $a,b \in \X$. Notice that even though the Bregman Lagrangian still involves the distance-generating function $h$, its natural motion is actually independent of $h$. Thus, the effect of $h$ is felt only via its interaction with $f$---this can also be seen in~\eqref{Eq:ELBreg2} where $h$ and $f$ only appear together in the final term. Furthermore, assuming $e^{\gamma_t} \to \infty$, the natural motion always converges to a limit point, which a priori can be anything. However, as we see from Theorem~\ref{Thm:Rate}, as soon as we introduce a convex potential function $f$, all motions converge to the minimizer $x^*$ of $f$.

\paragraph{Exponential convergence rate via uniform convexity}
\label{Sec:Exp}
In addition to the polynomial family in Section~\ref{Sec:Poly}, we can also study the subfamily of Bregman Lagrangians that have exponential convergence rates $O(e^{-ct})$, $c > 0$. As we discuss in Appendix~\ref{App:Exp}, in this case the link to discrete-time algorithms is not as clear. Using the same discretization technique as in Section~\ref{Sec:Poly} suggests that to get a matching convergence rate, constant progress is needed at each iteration.

From the discrete-time perspective, we show that the higher-order gradient algorithm~\eqref{Eq:AlgGrad} achieves an exponential convergence rate when the objective function $f$ is uniformly convex. Furthermore, we show that a restart scheme applied to the accelerated method~\eqref{Eq:AlgComp} achieves a better dependence on the condition number; this generalizes Nesterov's restart scheme for the case $p=3$~\cite[Section~5]{Nesterov08}. 

It is an open question to understand if there is a better connection between the discrete-time restart algorithms and the continuous-time exponential Lagrangian flows. In particular, it is of interest to consider whether a restart scheme is necessary to achieve exponential convergence in discrete time; we know it is not needed for the special case $p=2$, since a variant of Nesterov's accelerated gradient descent~\cite{Nesterov04} that incorporates the condition number also achieves the optimal convergence rate.

\section{Discussion}
\label{Sec:Discussion}

In this paper, we have presented a variational framework for understanding accelerated methods from a continuous-time perspective. 
We presented the general family of Bregman Lagrangian, which generates a family of second-order Lagrangian dynamics that minimize the objective function at an accelerated rate compared to gradient flows. These dynamics are related to each other by the operation of speeding up time, because the Bregman Lagrangian family is closed under time dilation. In the polynomial case, we showed how to discretize the second-order Lagrangian dynamics to obtain an accelerated algorithm with a matching convergence rate. The resulting algorithm accelerates a base algorithm by coupling it with a weighted mirror descent step. An example of a base algorithm is a higher-order gradient method, which in continuous time corresponds to a first-order rescaled gradient flow with a matching convergence rate.
Our continuous-time perspective makes clear that it is the mirror descent coupling that is more important for the acceleration phenomenon rather than the base algorithm. Indeed, the higher-order gradient algorithm operates on a smaller timescale than the enveloping mirror descent coupling step, so it makes no contribution in the continuous-time limit, and in principle we can use other base algorithms.

Our work raises many questions for further research. First, the case $p=2$ is worthy of further investigation. In particular, the assumptions needed to show convergence of the discrete-time algorithm ($\nabla^{p-1} f$ is Lipschitz) are different than those required to show existence and uniqueness of solutions of the continuous-time dynamics ($\nabla f$ is Lipschitz). In the case $p=2$ however, these assumptions match. This suggests a strong link between the discrete- and continuous-time dynamics that might help us understand why several results seem to be unique to the special case $p=2$.  Second, in discrete time, Nesterov's accelerated methods have been extended to various settings, for example to the stochastic setting. An immediate question is whether we can extend our Lagrangian framework to these settings. Third, we  would like to understand better the transition from continuous-time dynamics to discrete-time algorithms, and whether we can establish general assumptions that preserve desirable properties (e.g., convergence rate).  In Section~\ref{Sec:Poly} we saw that the polynomial convergence rate requires a higher-order smoothness assumption in discrete time, and in Section~\ref{Sec:FurtherProp} we discussed whether the exponential case requires a uniform convexity assumption. Finally, our work to date focuses on the convergence rates of the function values rather than the iterates. Recently there has been some work extending~\cite{SuBoydCandes14} to study the convergence of the iterates~\cite{AttouchPR15} and some perturbative aspects~\cite{AttouchC15}; it would be interesting to extend these results to the general Bregman Lagrangian.

At an abstract level, the general family of Bregman Lagrangian has a rich mathematical structure that deserves further study; we discussed some of these properties in Section~\ref{Sec:FurtherProp}. We hope that doing so will give us new insights into the nature of the optimization problem in continuous time, and help us design better dynamics with matching discrete-time algorithms.  For example, we can study how to use some of the appealing properties of the Hamiltonian formalism (e.g., volume preservation in phase space) to help us discretize the dynamics. We also wish to understand where the Bregman Lagrangian itself comes from, why it works so well,
and whether there are other Lagrangian families with similarly favorable properties.

\clearpage
\bibliography{draft_arxiv.bbl}

\clearpage
\appendix

\section{Proofs of results}

\subsection{Proof of Theorem~\ref{Thm:Time}}
\label{App:TimeProof}

The velocity and acceleration of the reparameterized curve $Y_t = X_{\tau(t)}$ are given by
\begin{align*}
\dot Y_t \,&=\, \dot \tau(t) \, \dot X_{\tau(t)} \\
\ddot Y_t \,&=\, \ddot \tau(t) \, \dot X_{\tau(t)} + \dot \tau(t)^2 \, \ddot X_{\tau(t)}.
\end{align*}
Inverting these relations, we get
\begin{subequations}\label{Eq:xyDot}
\begin{align}
\dot X_{\tau(t)} \,&=\, \frac{1}{\dot \tau(t)} \, \dot Y_t   \label{Eq:xyDot1} \\
\ddot X_{\tau(t)} \,&=\, \frac{1}{\dot \tau(t)^2} \ddot Y_t - \frac{\ddot \tau(t)}{\dot \tau(t)^3} \dot Y_t.   \label{Eq:xyDot2}
\end{align}
\end{subequations}

By assumption, the original curve $X_t$ satisfies the Euler-Lagrange equation~\eqref{Eq:ELBreg} for the Bregman Lagrangian $\L_{\alpha,\beta,\gamma}$. At time $\tau(t)$, this equation reads
\begin{equation*}
\begin{split} 
 &\ddot X_{\tau(t)} + (e^{\alpha_{\tau(t)}} - \dot \alpha_{\tau(t)}) \dot X_{\tau(t)} + e^{2\alpha_{\tau(t)} + \beta_{\tau(t)}}\Big[\nabla^2 h(X_{\tau(t)}+ e^{-\alpha_{\tau(t)}}\dot X_{\tau(t)})\Big]^{-1}\nabla f(X_{\tau(t)}) \\
& + e^{\alpha_{\tau(t)}}(\dot \gamma_{\tau(t)} - e^{\alpha_{\tau(t)}})\Big[\nabla^2 h(X_{\tau(t)} + e^{-\alpha_{\tau(t)}}\dot X_{\tau(t)})\Big]^{-1}(\nabla h(X_{\tau(t)} + e^{-\alpha_{\tau(t)}} \dot X_{\tau(t)}) - \nabla h(X_{\tau(t)})) = 0.
\end{split}
\end{equation*}
We now use the relations~\eqref{Eq:xyDot}. After multiplying by $\dot \tau(t)^2$ and collecting terms, we get
\begin{equation*}
\begin{split}
&\ddot Y_t + \left(\dot \tau(t) e^{\alpha_{\tau(t)}} - \dot \tau(t) \dot \alpha_{\tau(t)} - \frac{\ddot \tau(t)}{\dot \tau(t)}\right) \dot Y_t 
+ \dot \tau(t)^2 e^{2\alpha_{\tau(t)} + \beta_{\tau(t)}} \left[\nabla^2 h \left(Y_t + \frac{e^{-\alpha_{\tau(t)}}}{\dot \tau(t)} \dot Y_t\right)\right]^{-1} \nabla f(Y_t)  \\
&+ \dot \tau(t)^2 e^{\alpha_{\tau(t)}} (\dot \gamma_{\tau(t)} - e^{\alpha_{\tau(t)}}) \left[\nabla^2 h \left(Y_t + \frac{e^{-\alpha_{\tau(t)}}}{\dot \tau(t)} \dot Y_t\right)\right]^{-1} \left( \nabla h \left(Y_t + \frac{e^{-\alpha_{\tau(t)}}}{\dot \tau(t)} \dot Y_t\right) - \nabla h(Y_t)  \right) = 0.
\end{split}
\end{equation*}
Finally, with the definition of the modified parameters $\tilde\alpha,\tilde\beta,\tilde\gamma$~\eqref{Eq:ABC}, we can write this equation as
\begin{equation*}
\begin{split} 
 \ddot Y_t &+ (e^{\tilde\alpha_t} - \dot{\tilde\alpha}_t) \dot Y_t + e^{2\tilde\alpha_t + \tilde\beta_t}\Big[\nabla^2 h(Y_t + e^{-\tilde\alpha_t}\dot Y_t)\Big]^{-1}\nabla f(Y_t) \\
& + e^{\tilde\alpha_t}(\dot{\tilde\gamma}_t - e^{\tilde\alpha_t})\Big[\nabla^2 h(Y_t + e^{-\tilde\alpha_t}\dot Y_t)\Big]^{-1}(\nabla h(Y_t + e^{-\tilde\alpha_t} \dot Y_t) - \nabla h(Y_t)) = 0,
\end{split}
\end{equation*}
which we recognize as the Euler-Lagrange equation~\eqref{Eq:ELBreg} for the Bregman Lagrangian $\L_{\tilde\alpha,\tilde\beta,\tilde\gamma}$.

Furthermore, suppose $\alpha,\beta,\gamma$ satisfy the ideal scaling~\eqref{Eq:IdeSca}. Then
\begin{align*}
\dot{\tilde\beta}_t &= \frac{d}{dt} \beta_{\tau(t)} = \dot \tau(t) \dot \beta_{\tau(t)} \stackrel{\footnotesize\eqref{Eq:IdeScaBet}}{\le} \dot \tau(t) e^{\alpha_{\tau(t)}} = e^{\alpha_{\tau(t)} + \log \dot \tau(t)} = e^{\tilde \alpha_t} \\
\dot{\tilde\gamma}_t &= \frac{d}{dt} \gamma_{\tau(t)} = \dot \tau(t) \dot \gamma_{\tau(t)} \stackrel{\footnotesize\eqref{Eq:IdeScaGam}}{=} \dot \tau(t) e^{\alpha_{\tau(t)}} = e^{\alpha_{\tau(t)} + \log \dot \tau(t)} = e^{\tilde \alpha_t},
\end{align*}
which means that the modified parameters $\tilde\alpha,\tilde\beta,\tilde\gamma$ also satisfy the ideal scaling~\eqref{Eq:IdeSca}. The converse follows by considering the inverse function $\tau^{-1}(t)$ in place of $\tau(t)$.

$\hfill\qed$

\subsection{Existence and uniqueness of solution to the polynomial family}
\label{App:Existence}

In this section we discuss the existence and uniqueness of solution to the differential equation~\eqref{Eq:ELPoly} arising from the polynomial family of Bregman Lagrangian. 
We begin by writing the second-order equation~\eqref{Eq:ELPoly} as the pair of first-order equations~\eqref{Eq:ELPoly1}. We also write $W_t = \nabla h(Z_t)$, so we can write~\eqref{Eq:ELPoly1} as
\begin{subequations}\label{Eq:ELPoly2}
\begin{align}
\dot X_t &= \frac{p}{t} \left( \nabla h^*(W_t) - X_t \right)  \label{Eq:ELPoly2a} \\
\dot W_t &= -Cpt^{p-1}\nabla f(X_t)   \label{Eq:ELPoly2b}.
 \end{align}
\end{subequations}
Here $h^* \colon \X^* \to \R$ is the Legendre conjugate function of $h$, defined by
\begin{align}\label{Eq:hDual}
h^*(w) = \sup_{z \in \X} \left\{ \langle w,z \rangle - h(z) \right\},
\end{align}
where $\X^*$ is the dual space of $\X$, i.e., the space of all linear functionals over $\X$. Under the assumption that $h$ be essentially smooth, the supremum in~\eqref{Eq:hDual} is achieved by $z = \nabla h^*(w)$, and we have the relation that $\nabla h$ and $\nabla h^*$ are inverses of each other, i.e., $z = \nabla h^*(w) \Leftrightarrow w = \nabla h(z)$. Thus, with the definition $W_t = \nabla h(Z_t)$, we can write $Z_t = \nabla h^*(W_t)$, which gives us~\eqref{Eq:ELPoly2}.

Now assume $\nabla f$ and $\nabla h^*$ are Lipschitz continuous functions. Then over any bounded time intervals $[t_0, t_1]$ with $0 < t_0 < t_1$, the right-hand side of~\eqref{Eq:ELPoly2} is a Lipschitz continuous vector field. Thus, by the Cauchy-Lipschitz theorem, for any given initial conditions $(X_{t_0}, W_{t_0}) = (x_0, w_0)$ at time $t = t_0$, the system of differential equations~\eqref{Eq:ELPoly2} has a unique solution over the time interval $[t_0,t_1]$. 
Furthermore, the solution does not blow up in any finite time, since from Theorem~\ref{Thm:Rate} we know that the energy functional $\E_t$~\eqref{Eq:E} is non-increasing, so in particular, the Bregman divergence $D_h(x^*, X_t + \frac{t}{p} \dot X_t)$ is bounded above by a constant. 
Since $t_1$ is arbitrary, this shows that~\eqref{Eq:ELPoly2} has a unique maximal solution, i.e., $t_1$ can be extended to $t_1 \to +\infty$.

In the above argument we have started at time $t_0 > 0$, because the vector field in~\eqref{Eq:ELPoly2} has a singularity at $t = 0$. 
For $p=2$, Su et al.~\cite{SuBoydCandes14} and Krichene et al.~\cite{Krichene15} treat the case when we start at $t=0$ with initial condition $(X_0, W_0) = (x_0, \nabla h(x_0))$, so that $\dot X_0 = 0$. In that case, they show that the system~\eqref{Eq:ELPoly2} still has a unique solution for all time $[0,\infty)$, by replacing the $p/t$ coefficient by the approximation $p/\max\{t,\delta\}$ for $\delta > 0$ and letting $\delta \to 0$. We can adapt this technique to the more general case~\eqref{Eq:ELPoly2}; alternatively, we can appeal to the time dilation property and state that since the general system~\eqref{Eq:ELPoly2} is the result of speeding up the $p=2$ case by time dilation function $\tau(t) = t^{p/2}$, once we know a unique solution exists for $p=2$, we can also conclude that it exists for all $p > 0$.

\subsection{Proof of Theorem~\ref{Thm:AlgRate}}
\label{Thm:AlgRateProof}

We define the following function, which is a generalization of Nesterov's {\em estimate function} from~\cite{Nesterov08},
\begin{align}\label{Eq:EstFunc}
\psi_k(x) = Cp \sum_{i=0}^k i^{(p-1)} \big[f(y_i) + \langle \nabla f(y_i), x-y_i \rangle \big] + \frac{1}{\epsilon} D_h(x,x_0).
\end{align}
The estimate function $\psi_k$ arises as the objective function that the sequence $z_k$ is optimizing in~\eqref{Eq:Algz}. Indeed, the optimality condition for the $z_k$ update~\eqref{Eq:Algz} is
\begin{align*}
\nabla h(z_k) = \nabla h(z_{k-1}) - \epsilon Cpk^{(p-1)} \nabla f(y_k).
\end{align*}
By unrolling the recursion, we can write
\begin{align*}
\nabla h(z_k) = \nabla h(z_0) - \epsilon Cp \sum_{i=0}^k i^{(p-1)} \nabla f(y_i),
\end{align*}
and since $x_0 = z_0$, we can write this equation as $\nabla \psi_k(z_k) = 0$. Since $\psi_k$ is a convex function, this means $z_k$ is the minimizer of $\psi_k$. Thus, we can equivalently write the update for $z_k$ as
\begin{align}\label{Eq:EstFunc0}
z_k = \arg\min_z \: \psi_k(z).
\end{align}

For proving the convergence rate for the algorithm~\eqref{Eq:Alg}, we have the following property.
\begin{lemma}\label{Lem:EstFunc}
For all $k \ge 0$, we have
\begin{align}\label{Eq:EstFunc1}
\psi_k(z_k) \,\ge\, Ck^{(p)} f(y_k).
\end{align}
\end{lemma}
\begin{proof}
We proceed via induction on $k \ge 0$. The base case $k = 0$ is true since both sides equal zero. Now assume~\eqref{Eq:EstFunc1} holds for some $k \ge 0$; we will show it also holds for $k+1$.

Since $h$ is $1$-uniformly convex of order $p$, the rescaled Bregman divergence $\frac{1}{\epsilon} D_h(x,x_0)$ is $(\frac{1}{\epsilon})$-uniformly convex. Thus, the estimate function $\psi_k$~\eqref{Eq:EstFunc} is also $(\frac{1}{\epsilon})$-uniformly convex of order $p$.
Since $z_k$ is the minimizer of $\psi_k$, $\nabla \psi_k(z_k) = 0$, so for all $x \in \X$ we have
\begin{align*}
\psi_k(x) \,=\, \psi_k(z_k) + D_{\psi_k}(x,z_k) \,\ge\, \psi_k(z_k) + \frac{1}{\epsilon p} \|x-z_k\|^p.
\end{align*}
Applying the inductive hypothesis~\eqref{Eq:EstFunc1} and using the convexity of $f$ gives us
\begin{align*}
\psi_k(x) 
\,\ge\, Ck^{(p)} \big[ f(y_{k+1}) + \langle \nabla f(y_{k+1}), y_k-y_{k+1} \rangle \big] + \frac{1}{\epsilon p} \|x-z_k\|^p.
\end{align*}
We now add $Cp(k+1)^{(p-1)}[f(y_{k+1}) + \langle \nabla f(y_{k+1}), x-y_{k+1}\rangle]$ to both sides of the equation to obtain
\begin{align}\label{Eq:EstFunc3}
\psi_{k+1}(x)
\,\ge\,
C(k+1)^{(p)} \big[ f(y_{k+1}) + \big\langle \nabla f(y_{k+1}),\, x_{k+1}-y_{k+1} + \tau_k (x-z_k) \big\rangle \big] + \frac{1}{\epsilon p} \|x-z_k\|^p,
\end{align}
where $\tau_k = \frac{p (k+1)^{(p-1)}}{(k+1)^{(p)}} = \frac{p}{k+p}$, and where we have also used the definition of $x_{k+1}$ as a convex combination of $y_k$ and $z_k$ with weight $\tau_k$~\eqref{Eq:Algx}.

Note that the first term in~\eqref{Eq:EstFunc3} gives our desired inequality~\eqref{Eq:EstFunc1} for $k+1$. So to finish the proof, we have to prove the remaining terms in~\eqref{Eq:EstFunc3} are nonnegative.
We do so by applying two inequalities. We first apply the inequality~\eqref{Eq:yIneq} to the term $\langle \nabla f(y_{k+1}), x_{k+1}-y_{k+1} \rangle$, so from~\eqref{Eq:EstFunc3} we have
\begin{equation}\label{Eq:EstFunc4}
\begin{split}
\psi_{k+1}(x)
\,&\ge\,
C(k+1)^{(p)} f(y_{k+1}) \,+\, C(k+1)^{(p)} M \epsilon^{\frac{1}{p-1}} \, \|\nabla f(y_{k+1})\|_*^{\frac{p}{p-1}} \\
&\qquad\qquad+\, Cp(k+1)^{(p-1)} \langle \nabla f(y_{k+1}), x-z_k \rangle \,+\, \frac{1}{\epsilon p} \|x-z_k\|^p.
\end{split}
\end{equation}
Next, we apply the Fenchel-Young inequality~\cite[Lemma~2]{Nesterov08}
\begin{align}\label{Eq:FenchelYoung}
\langle s,u \rangle + \frac{1}{p} \|u\|^p \ge - \frac{p-1}{p} \|s\|^{\frac{p}{p-1}}_*
\end{align}
with the choices $u = \epsilon^{-\frac{1}{p}} (x-z_k)$ and $s =  \epsilon^{\frac{1}{p}} Cp  (k+1)^{(p-1)} \nabla f(y_{k+1})$. Then from~\eqref{Eq:EstFunc4}, we obtain
\begin{align}
\psi_{k+1}(x)
\,\ge\,
C(k+1)^{(p)}  \left[ f(y_{k+1})
+ \left( M - \frac{p-1}{p} \, p^{\frac{p}{p-1}}\, C^{\frac{1}{p-1}} \frac{\{(k+1)^{(p-1)}\}^{\frac{p}{p-1}}}{(k+1)^{(p)}} \right) \, \epsilon^{\frac{1}{p-1}} \, \|\nabla f(y_{k+1})\|_*^{\frac{p}{p-1}} \right]. \notag
\end{align}
Notice that $\{(k+1)^{(p-1)}\}^{\frac{p}{p-1}} \le (k+1)^{(p)}$. Then from the assumption $C \le M^{p-1}/p^p$, we see that the second term inside the parentheses  is nonnegative. Hence we conclude the desired inequality $\psi_{k+1}(x) \ge C (k+1)^{(p)} f(y_{k+1})$. Since $x \in \X$ is arbitrary, it also holds for the minimizer $x = z_{k+1}$ of $\psi_{k+1}$, finishing the induction. 
\end{proof}

With Lemma~\ref{Lem:EstFunc} in hand, we can complete the proof of Theorem~\ref{Thm:AlgRate}.

\begin{proof}[Proof of Theorem~\ref{Thm:AlgRate}]
Since $f$ is convex, we can bound the estimate sequence $\psi_k$ by
\begin{align*}
\psi_k(x) \,\le\, Cp \sum_{i=0}^k i^{(p-1)} f(x) + \frac{1}{\epsilon} D_h(x,x_0) \,=\, Ck^{(p)} f(x) + \frac{1}{\epsilon} D_h(x,x_0).
\end{align*}
This holds for all $x \in \X$, and in particular for the minimizer $x^*$ of $f$. Combining the bound with the result of Lemma~\ref{Lem:EstFunc}, and recalling that $z_k$ is the minimizer of $\psi_k$, we get
\begin{align*}
Ck^{(p)}f(y_k) \,\le\, \psi_k(z_k) \,\le\, \psi_k(x^*) \,\le\, Ck^{(p)} f(x^*) + \frac{1}{\epsilon} D_h(x^*,x_0).
\end{align*}
Rearranging and dividing by $Ck^{(p)}$ gives us the desired convergence rate~\eqref{Eq:AlgRate}.
\end{proof}

\subsection{Proof of Lemma~\ref{Lem:Update}}
\label{App:UpdateProof}

We follow the approach of~\cite[Lemma~6]{Nesterov08}. Since $y$ solves the optimization problem~\eqref{Eq:Update}, it satisfies the optimality condition
\begin{align}\label{Eq:UpdateProof1}
\sum_{i=1}^{p-1} \frac{1}{(i-1)!} \nabla^i f(x) \, (y-x)^{i-1} + \frac{N}{\epsilon} \|y-x\|^{p-2} \, (y-x) = 0.
\end{align}
Furthermore, since $\nabla^{p-1} f$ is $\frac{(p-1)!}{\epsilon}$-Lipschitz, we have the following error bound on the $(p-2)$-nd order Taylor expansion of $\nabla f$,
\begin{align}\label{Eq:UpdateProof2}
\left\|\nabla f(y) - \sum_{i=1}^{p-1} \frac{1}{(i-1)!} \nabla^i f(x) \, (y-x)^{i-1}\right\|_* \le \frac{1}{\epsilon} \|y-x\|^{p-1}.
\end{align}
Substituting~\eqref{Eq:UpdateProof1} to~\eqref{Eq:UpdateProof2} and writing $r = \|y-x\|$, we obtain
\begin{align}\label{Eq:UpdateProof2a}
\left\|\nabla f(y) + \frac{Nr^{p-2}}{\epsilon} \, (y-x)\right\|_* \,\le\, \frac{r^{p-1}}{\epsilon}.
\end{align}
Squaring both sides, expanding, and rearranging the terms, we get the inequality
\begin{align}\label{Eq:UpdateProof3}
\langle \nabla f(y), x-y \rangle
\,\ge\, \frac{\epsilon}{2Nr^{p-2}} \|\nabla f(y)\|_*^2 + \frac{(N^2-1)r^p}{2N\epsilon}.
\end{align}
Note that if $p=2$, then the first term in~\eqref{Eq:UpdateProof3} already implies the desired bound~\eqref{Eq:UpdateIneq}. Now assume $p \ge 3$. The right-hand side of~\eqref{Eq:UpdateProof3} is of the form $A/r^{p-2} + Br^p$, which is a convex function of $r > 0$ and minimized by $r^* = \left\{\frac{(p-2)}{p} \frac{A}{B} \right\}^{\frac{1}{2p-2}}$, yielding a minimum value of
\begin{align*}
\frac{A}{(r^*)^{p-2}} + B(r^*)^p
\,=\, A^{\frac{p}{2p-2}} B^{\frac{p-2}{2p-2}} \left[\left(\frac{p}{p-2}\right)^{\frac{p-2}{2p-2}} + \left(\frac{p-2}{p}\right)^{\frac{p}{p-2}}\right]
\,\ge\, A^{\frac{p}{2p-2}} B^{\frac{p-2}{2p-2}}.
\end{align*}
Substituting the values $A = \frac{\epsilon}{2N} \|\nabla f(y)\|_*^2$ and $B = \frac{1}{2N\epsilon} (N^2-1)$ from~\eqref{Eq:UpdateProof3}, we obtain
\begin{align*}
\langle \nabla f(y), x-y \rangle
\,&\ge\, \left(\frac{\epsilon}{2N} \|\nabla f(y)\|_*^2\right)^{\frac{p}{2p-2}} \left(\frac{1}{2N\epsilon} (N^2-1)\right)^{\frac{p-2}{2p-2}}
\,=\, \frac{(N^2-1)^{\frac{p-2}{2p-2}}}{2N} \epsilon^{\frac{1}{p-1}} \|\nabla f(y)\|_*^{\frac{p}{p-1}}
\end{align*}
which proves~\eqref{Eq:UpdateIneq}.

To obtain the first inequality of~\eqref{Eq:UpdateIneq2}, we use Cauchy-Schwarz inequality on~\eqref{Eq:UpdateIneq},
\begin{align*}
\frac{(N^2-1)^{\frac{p-2}{2p-2}}}{2N} \epsilon^{\frac{1}{p-1}} \|\nabla f(y)\|_*^{\frac{p}{p-1}}
\,\le\, \langle \nabla f(y), x-y \rangle
\,\le\, \|\nabla f(y)\|_* \, \|x-y\|
\end{align*}
and cancel out $\|\nabla f(y)\|_*$ from both sides.
For the second inequality of~\eqref{Eq:UpdateIneq2}, we use triangle inequality on the left hand side of~\eqref{Eq:UpdateProof2a},
\begin{align*}
\frac{N r^{p-1}}{\epsilon} - \|\nabla f(y)\|_* \,\le\, \left\|\nabla f(y) + \frac{Nr^{p-2}}{\epsilon} \, (y-x)\right\|_* \,\le\, \frac{r^{p-1}}{\epsilon}.
\end{align*}
Rearranging the terms and taking the $(p-1)$-st root of both sides gives us the result~\eqref{Eq:UpdateIneq2}.
\qed

\subsection{Proof of Theorem~\ref{Thm:AlgGradRate}}
\label{App:AlgGradRateProof}

This proof follows the approach in the proof of~\cite[Theorem~1]{Nesterov08}.
We first prove the following lemma. Here $\delta_k = f(x_k) - f(x^*) \ge 0$ denotes the residual value at iteration $k$.

\begin{lemma}\label{Lem:AlgGradRate}
Under the setting of Theorem~\ref{Thm:AlgGradRate}, we have
\begin{align}\label{Eq:AlgGradResidual}
\delta_{k+1} \le \delta_k - \frac{(p-1)}{p} \cdot \left(\frac{\epsilon \delta_k^p}{(N+1) R^p}\right)^{\frac{1}{p-1}}. 
\end{align}
\end{lemma}
\begin{proof}
Since $f$ is $\frac{(p-1)!}{\epsilon}$-smooth of order $p-1$, by the Taylor remainder theorem we have the bound
\begin{align*}
|f_{p-1}(x;x_k) - f(x)| \le \frac{1}{\epsilon p} \|x-x_k\|^p.
\end{align*}
Then from the definition of $x_{k+1}$~\eqref{Eq:AlgGrad}, we have
\begin{align}\label{Eq:AlgGradResidual2}
f(x_{k+1}) \,&=\, \min_{x \in \X} \left\{ f_{p-1}(x;x_k) + \frac{N}{\epsilon p} \|x-x_k\|^p \right\} 
\le\, \min_{x \in \X} \left\{ f(x) + \frac{N+1}{\epsilon p} \|x-x_k\|^p \right\}.
\end{align}
Plugging in $x = x_k$ on the right-hand side of~\eqref{Eq:AlgGradResidual2} shows that $f(x_{k+1}) \le f(x_k)$; that is, the algorithm~\eqref{Eq:AlgGrad} is a descent method. In particular, for all $k \ge 0$ we have $\|x_k - x^*\| \le R$, where $R = \sup_{x \colon f(x) \le f(x_0)} \|x-x^*\|$ is the radius of the level set as defined in Theorem~\ref{Thm:AlgGradRate}.
Moreover, plugging in $x = x^*$ on the right-hand side of~\eqref{Eq:AlgGradResidual2} gives us
\begin{align}\label{Eq:AlgGradResidual3}
f(x_{k+1}) - f(x^*) \le \frac{N+1}{\epsilon p} \|x_k - x^*\|^p \le \frac{N+1}{\epsilon p} R^p.
\end{align}

Now for any $\lambda \in [0,1]$, consider the midpoint
\begin{align*}
x_\lambda = x^* + (1-\lambda)(x_k-x^*) = \lambda x^* + (1-\lambda) x_k.
\end{align*}
By Jensen's inequality, $f(x_\lambda) \le \lambda f(x^*) + (1-\lambda) f(x_k)$. We also have $\|x_\lambda - x_k\| = \lambda \|x_k - x^*\| \le \lambda R$. Plugging in the point $x_\lambda$ to the right-hand side of~\eqref{Eq:AlgGradResidual2} gives
\begin{align*}
f(x_{k+1})
\,\le\, f(x_\lambda) + \frac{N+1}{\epsilon p} \|x_\lambda-x_k\|^p
\,\le\, \lambda f(x^*) + (1-\lambda) f(x_k) + \frac{N+1}{\epsilon p} R^p \lambda^p. 
\end{align*}
With the notation $\delta_k = f(x_k) - f(x^*)$, we can write the last inequality as
\begin{align}\label{Eq:AlgGradResidual4}
\delta_{k+1} \:\le\: (1-\lambda) \delta_k + \frac{N+1}{\epsilon p} R^p \lambda^p.
\end{align}
The right-hand side is a convex function of $\lambda$, which is minimized at $\lambda^* = \left\{ \frac{\epsilon}{N+1} \frac{\delta_k}{R^p} \right\}^{\frac{1}{p-1}}$. Note that $\lambda^* \in [0,1]$ by~\eqref{Eq:AlgGradResidual3}. Plugging in $\lambda^*$ to~\eqref{Eq:AlgGradResidual4} yields the desired bound~\eqref{Eq:AlgGradResidual}.
\end{proof}

With Lemma~\ref{Lem:AlgGradRate}, we can complete the proof of Theorem~\ref{Thm:AlgGradRate}.

\begin{proof}[Proof of Theorem~\ref{Thm:AlgGradRate}]
Define the energy functional $e_k = \delta_k^{-\frac{1}{p-1}}$. We can write
\begin{align}\label{Eq:AlgGradProof}
e_{k+1} - e_k
  \,=\, \frac{1}{ \delta_{k+1}^{\frac{1}{p-1}}} - \frac{1}{ \delta_k^{\frac{1}{p-1}}} 
  \,=\,  \frac{\delta_{k}^{\frac{1}{p-1}}-\delta_{k+1}^{\frac{1}{p-1}}}{\delta_{k+1}^{\frac{1}{p-1}} \cdot \delta_{k}^{\frac{1}{p-1}}}
  \,=\, \frac{\delta_k-\delta_{k+1}}{\delta_{k+1}^{\frac{1}{p-1}} \cdot \delta_{k}^{\frac{1}{p-1}}} \cdot \frac{1}{\left( \sum_{i=0}^{p-2} \delta_k^{\frac{i}{p-1}} \cdot \delta_{k+1}^{\frac{p-2-i}{p-1}} \right)}.
\end{align}
Since $\delta_{k+1} \le \delta_k$, we can upper bound the summation in the denominator of~\eqref{Eq:AlgGradProof} by $(p-1) \delta_k^{\frac{p-2}{p-1}}$. We use Lemma~\ref{Lem:AlgGradRate} to lower bound~$\delta_k-\delta_{k-1}$, obtaining
\begin{align}\label{Eq:AlgGradProof2}
e_{k+1}-e_k
  \:\ge\: \frac{(p-1)}{p} \cdot \left(\frac{\epsilon \delta_k^{p}}{(N+1)R^{p}}\right)^{\frac{1}{p-1}} \cdot \frac{1}{\delta_k^{\frac{2}{p-1}}} \cdot \frac{1}{(p-1) \delta_k^{\frac{p-2}{p-1}}}
  \:=\: \frac{1}{p} \cdot \left(\frac{\epsilon}{(N+1)R^p}\right)^{\frac{1}{p-1}}.
\end{align}
Summing~\eqref{Eq:AlgGradProof2} and telescoping the terms, we get
\begin{align*}
\frac{1}{(f(x_k) - f(x^*))^{\frac{1}{p-1}}} \,=\, e_k \,\ge\, e_k-e_0 \,\ge\, \frac{k}{p} \cdot \left(\frac{\epsilon}{(N+1)R^p}\right)^{\frac{1}{p-1}}
\end{align*}
which gives us the desired conclusion~\eqref{Eq:AlgGradRate}. 
\end{proof}

\subsection{Proof of Theorem~\ref{Thm:RescGradFlow}}
\label{App:RescGradFlowProof}

We write the higher-order gradient algorithm~\eqref{Eq:AlgGrad} (with $N=1$) as
    \begin{align}\label{Eq:AlgGrad2}
    x_{k+1} - x_k = \arg\min_u \left\{ f(x_k) + \langle \nabla f(x_k), u \rangle + \cdots + \frac{1}{(p-1)!} \nabla^{p-1} f(x_k) u^{p-1} + \frac{1}{\epsilon p} \|u\|^p \right\}.
    \end{align}
Our goal is to express the sequence $x_k$ as a discretization $x_k = X_t$, $x_{k+1} = X_{t+\delta} \approx X_t + \delta \dot X_t$ of some continuous-time curve $X_t$ with time step $\delta > 0$, which will be a function of $\epsilon$. To that end, we write $u = \delta v$, so~\eqref{Eq:AlgGrad2} becomes
    \begin{align*}
    \frac{x_{k+1} - x_k}{\delta} = \arg\min_v \left\{ f(x_k) + \delta \langle \nabla f(x_k), v \rangle + \cdots + \frac{\delta^{p-1}}{(p-1)!} \nabla^{p-1} f(x_k) v^{p-1} + \frac{\delta^p}{\epsilon p} \|v\|^p \right\}.
    \end{align*}
Eliminating the constant term $f(x_k)$ from the right-hand side, which does not change the minimizer, and canceling a factor of $\delta$, we get
    \begin{align*}
    \frac{x_{k+1} - x_k}{\delta} = \arg\min_v \left\{ \langle \nabla f(x_k), v \rangle + \frac{\delta}{2} \nabla^2 f(x_k) v^2+ \cdots + \frac{\delta^{p-2}}{(p-1)!} \nabla^{p-1} f(x_k) v^{p-1} + \frac{\delta^{p-1}}{\epsilon p} \|v\|^p \right\}.
    \end{align*}
We see that the first term in the objective function does not depend on $\delta$. As $\epsilon \to 0$, for the equation to have a meaningful limit, we have to set $\delta^{p-1} = \epsilon$, so the last term in the objective function becomes a constant. On the other hand, the middle terms all have dependence on $\delta = \epsilon^{\frac{1}{p-1}}$, so as $\epsilon \to 0$, those terms vanish. Thus, the limit as $\epsilon \to 0$ is
    \begin{align}\label{Eq:AlgGrad3}
    \dot X_t = \arg\min_v \left\{ \langle \nabla f(X_t), v \rangle + \frac{1}{p} \|v\|^p \right\}.
    \end{align}
Equivalently, $\dot X_t$ satisfies the optimality condition
    \begin{align}\label{Eq:AlgGrad4}
    \nabla f(X_t) + \|\dot X_t\|^{p-2} \dot X_t = 0.
    \end{align}
This gives us the relation $\|\nabla f(X_t)\|_* = \|\dot X_t\|^{p-1}$, so we can also write~\eqref{Eq:AlgGrad4} as
    \begin{align*}
    \dot X_t = -\frac{\nabla f(X_t)}{\|\dot X_t\|^{p-2}} = -\frac{\nabla f(X_t)}{\|\nabla f(X_t)\|_*^{\frac{p-2}{p-1}}},
    \end{align*}
which is the rescaled gradient flow as claimed in~\eqref{Eq:RescGradFlow}.

We note that the rescaled gradient flow~\eqref{Eq:RescGradFlow} is a descent method, since
    \begin{align*}
    \frac{d}{dt} f(X_t) = \langle \nabla f(X_t), \dot X_t \rangle = -\|\nabla f(X_t)\|_*^{\frac{p}{p-1}} \le 0.
    \end{align*}
Now to establish the convergence rate of the rescaled gradient flow~\eqref{Eq:RescGradFlow}, we consider the energy functional
    \begin{align}\label{Eq:RescGradFlowE}
    \E_t = (f(X_t) - f(x^*))^{-\frac{1}{p-1}}
    \end{align}
which is the same energy functional as in the discrete-time convergence proof in Appendix~\ref{App:AlgGradRateProof}. The energy functional $\E_t$ has time derivative
    \begin{align*}
    \dot \E_t
     = -\frac{1}{(p-1)} \frac{\langle \nabla f(X_t), \dot X_t \rangle}{(f(X_t)-f(x^*))^{\frac{p}{p-1}}}.
    \end{align*}
If $X_t$ satisfies the rescaled gradient flow equation~\eqref{Eq:RescGradFlow}, then $\dot \E_t$ simplifies to
    \begin{align}\label{Eq:RescGradFlowEdot}
    \dot \E_t
      = \frac{1}{(p-1)} \left(\frac{\|\nabla f(X_t)\|_*}{f(X_t)-f(x^*)}\right)^{\frac{p}{p-1}}.
    \end{align}
By the convexity of $f$ and the Cauchy-Schwarz inequality, we have
    \begin{align*}
    0 \le f(X_t)-f(x^*) \le \langle \nabla f(X_t), X_t-x^* \rangle \le \|\nabla f(X_t)\|_* \, \|X_t-x^*\|.
    \end{align*}
Since the rescaled gradient flow is a descent method, we have $\|X_t - x^*\| \le R$. Therefore, from~\eqref{Eq:RescGradFlowEdot} we get the bound
    \begin{align*}
    \dot \E_t
      \geq \frac{1}{(p-1)} \frac{1}{\|X_t-x^*\|^{\frac{p}{p-1}}}
      \geq \frac{1}{(p-1) R^{\frac{p}{p-1}}}.
    \end{align*}
This means that $\E_t$ increases at least linearly, so
    \begin{align*}
    \frac{1}{(f(X_t) - f(x^*))^{\frac{1}{p-1}}} \,=\, \E_t \,\geq\, \E_0 + \frac{t}{(p-1) \, R^{\frac{p}{p-1}}} \,\ge\, \frac{t}{(p-1) \, R^{\frac{p}{p-1}}},
    \end{align*}
which gives us the desired result~\eqref{Eq:RescGradFlowRate}.
\qed

\paragraph{Remark:}
From the proof above, we see that rescaled gradient flow~\eqref{Eq:RescGradFlow} is a generalization of the usual gradient flow (the case $p=2$) which is obtained by replacing the squared norm by the $p$-th power of the norm in the variational formulation~\eqref{Eq:AlgGrad3}.
It turns out that when the objective function is the $p$-th power of the norm, $f(x) = \frac{1}{p} \|x\|^p$, the rescaled gradient flow~\eqref{Eq:RescGradFlow} reduces to an explicit equation. Specifically, in this case we have $\nabla f(x) = \|x\|^{p-2} x$, so $\|\nabla f(x)\|_* = \|x\|^{p-1}$. Therefore, the rescaled gradient flow equation~\eqref{Eq:RescGradFlow} becomes
\begin{align*}
\dot X_t = -\frac{\nabla f(X_t)}{\|\nabla f(X_t)\|_*^{\frac{p-2}{p-1}}} = -\frac{\|X_t\|^{p-2} X_t}{\|X_t\|^{p-2}} = -X_t,
\end{align*}
which is now independent of $p$, and has an explicit solution $X_t = e^{-t} X_0$.

\paragraph{Alternative proof of convergence rate.}
In the proof above, we can also use the following alternative energy functional,
\begin{align}\label{Eq:EAlt}
\tilde{\E_t} = t^p(f(X_t) - f(x^*)).
\end{align}
Its time derivative is
\begin{subequations}
\begin{align}
\dot{\tilde{\E}}_t \,&=\, pt^{p-1}(f(X_t) - f(x^*)) + t^p \langle \nabla f(X_t), \dot X_t \rangle \notag \\
&\le\, pt^{p-1} \langle \nabla f(X_t), X_t - x^* \rangle + t^p \langle \nabla f(X_t), \dot X_t \rangle \label{Eq:EAltdot1} \\
&=\, pt^{p-1} \langle \nabla f(X_t), X_t - x^* \rangle - t^p \| \nabla f(X_t) \|^{\frac{p}{p-1}}_*,
\label{Eq:EAltdot2}
\end{align}
\end{subequations}
where~\eqref{Eq:EAltdot1} follows from the convexity of $f$, and in~\eqref{Eq:EAltdot2} we have substituted the rescaled gradient flow dynamic~\eqref{Eq:RescGradFlow}.
We now apply the Fenchel-Young inequality~\eqref{Eq:FenchelYoung} with $s = t^{p-1} \nabla f(X_t)$ and $u = -(p-1) (X_t-x^*)$, to obtain
\begin{align}\label{Eq:EAltdot3}
\dot{\tilde{\E}}_t \,\le\, \frac{1}{p-1} \|(p-1)(X_t - x^*)\|^p \,\le\, (p-1)^{p-1} R^p,
\end{align}
where in the last step we have used the fact that $\|X_t - x^*\| \le R$ since rescaled gradient flow is a descent method. Integrating~\eqref{Eq:EAltdot3} and plugging in the definition of $\tilde{\E_t}$~\eqref{Eq:EAlt}, we obtain
\begin{align*}
f(X_t) - f(x^*) \,\le\, \frac{(p-1)^{p-1} \, R^p}{t^{p-1}},
\end{align*}
which is exactly the same bound as claimed in~\eqref{Eq:RescGradFlowRate}.

\section{Further properties}

\subsection{Exponential convergence rate via uniform convexity}
\label{App:Exp}

Similar to the polynomial case in Section~\ref{Sec:Poly}, in this section we study the subfamily of Bregman Lagrangian~\eqref{Eq:BregLag} with the following choice of parameters, parameterized by $c>0$,
\begin{subequations}\label{Eq:ABCExp}
\begin{align}
\alpha_t &=  \log c\\
\beta_t &= ct \\
\gamma_t &= ct.
\end{align}
\end{subequations}
The parameters~\eqref{Eq:ABCExp} satisfy the ideal scaling condition~\eqref{Eq:IdeSca}, with an equality on the first condition~\eqref{Eq:IdeScaBet}. The Euler-Lagrange equation~\eqref{Eq:ELBreg2} in this case is given by
\begin{align}\label{Eq:ELExp}
\ddot X_t + c \dot X_t + c^2 e^{ct} \left[\nabla^2 h\left(X_t + \frac{1}{c} \dot X_t\right)\right]^{-1} \nabla f(X_t) = 0,
\end{align}
and by Theorem~\ref{Thm:Rate}, it has an $O(e^{-ct})$ rate of convergence. 
Thus, whereas the polynomial Lagrangian flow~\eqref{Eq:ELPoly} has a polynomial rate of convergence, the exponential Lagrangian flow~\eqref{Eq:ELExp} has an exponential rate of convergence.
Furthermore, from the time-dilation property in Theorem~\ref{Thm:Time}, we see that we can obtain the exponential curve~\eqref{Eq:ELExp} by speeding up the polynomial curve~\eqref{Eq:ELPoly} using a time-dilation function $\tau(t) = e^{ct/p}$.

However, unlike the polynomial Lagrangian flow~\eqref{Eq:ELPoly}, the process of discretizing the exponential Lagrangian flow~\eqref{Eq:ELExp} is not as straightforward. Following the same approach as the polynomial family, we write the second-order equation \eqref{Eq:ELExp} as the following pair of first-order equations: 
\begin{subequations}\label{Eq:ELExp1}
 \begin{align}
 Z_t &= X_t + \frac{1}{c}\dot X_t    \label{Eq:ELExp1a} \\
\frac{d}{dt} \nabla h(Z_t) &= -ce^{ct}\nabla f(X_t)   \label{Eq:ELExp1b}.
 \end{align}
\end{subequations}
Now we discretize $X_t$ and $Z_t$ into sequences $x_k$ and $z_k$ with time step $\delta >0$, so that $t = \delta k$ as before. In doing so, we can write~\eqref{Eq:ELExp1} as the following discrete-time equations similar to \eqref{Eq:Alg0x} and \eqref{Eq:Alg0z}:
\begin{subequations}\label{Eq:expFlowAlg}
\begin{align}
x_{k+1} &= c\delta z_k + (1 - c\delta) x_k \label{Eq:expFlowAlga} \\
z_{k+1} &= \arg\min_z \,\left\{ ce^{c\delta k} \langle \nabla f(x_k), z \rangle + \frac{1}{\delta} D_h(z,z_k) \right\} \label{Eq:expFlowAlgb}.
\end{align}
\end{subequations}
Note that the weight in~\eqref{Eq:expFlowAlga} is independent of time, but depends on $\delta$, and~\eqref{Eq:expFlowAlgb} suggests the step size $\epsilon = \delta$ in the algorithm.
If our analogy between continuous and discrete-time convergence holds, then given the $O(e^{-ct})$ convergence rate in continuous time, we expect a matching $O(\frac{1}{\epsilon} e^{-ck})$ convergence rate in discrete time. However, it is not clear how to obtain that rate via~\eqref{Eq:expFlowAlg}. If we try to adapt the proof of Theorem~\ref{Thm:AlgRate}, we find that in order to conclude a convergence rate $O(\delta e^{-c\delta k})$, we need to introduce a sequence $y_k$ satisfying the following analog of inequality~\eqref{Eq:yIneq} (with the ideal choice $p=\infty$):
\begin{align}\label{Eq:expNeeded}
\langle \nabla f(y_k), x_k-y_k \rangle \,\ge\, M \|\nabla f(y_k)\|_*.
\end{align}
Notice that the rates are consistent if we set $\epsilon = \delta = 1$. However, the condition~\eqref{Eq:expNeeded} means we need to make a constant improvement in each iteration from $x_k$ to $y_k$, although we are also free on how we choose to construct $y_k$ and impose any assumptions on $f$.

In the remainder of this section, we approach this problem from a discrete-time perspective, and study the performance of the higher-order gradient algorithm~\eqref{Eq:AlgGrad} and its accelerated variant~\eqref{Eq:AlgComp} when $f$ is uniformly convex.

\subsubsection{Exponential convergence rate of higher-order gradient algorithm}
\label{App:ExpConvergenceGrad}

In this section we show that the higher-order gradient algorithm~\eqref{Eq:AlgGrad} has an exponential convergence rate when the objective function $f$ is uniformly convex of order $p \ge 2$; this generalizes the results in~\cite[Section~5]{Nesterov08} for the case $p=3$, and the classical result of gradient descent for the case $p=2$~\cite{Nesterov04}.

Specifically, we have the following result. Recall the definition of smoothness in~\eqref{Eq:Smooth}, and the definition of uniform convexity in~\eqref{Eq:UnifConv}.

\begin{theorem}\label{Thm:GradRateUnif}
Suppose $f$ is $\frac{(p-1)!}{\epsilon}$-smooth of order $p-1$, and $\sigma$-uniformly convex of order $p$. Then the $p$-th order gradient algorithm~\eqref{Eq:AlgGrad} with $N > 1$ has convergence rate
\begin{align}\label{Eq:GradRateUnif}
f(x_{k+1}) - f(x^*) \le \frac{(N+1)\|x_0-x^*\|^p}{\epsilon p \big(1+L\kappa^{\frac{1}{p-1}}\big)^k}
\,=\, O\left(\frac{1}{\epsilon}  \, \exp(-L\kappa^{\frac{1}{p-1}} k)\right),
\end{align}
where $L = (N^2-1)^{\frac{p-2}{2p-2}}/(2N)$, and $\kappa = \epsilon \sigma$ is the inverse condition number (which we assume is small).
\end{theorem}
\begin{proof}
By inequality~\eqref{Eq:UpdateIneq} from Lemma~\ref{Lem:Update}, we know that since $f$ is $\frac{(p-1)!}{\epsilon}$-smooth of order $p-1$, 
\begin{align*}
\langle \nabla f(x_{k+1}), x_k-x_{k+1} \rangle \,\ge\, L \epsilon^{\frac{1}{p-1}} \|\nabla f(x_{k+1})\|_*^{\frac{p}{p-1}},
\end{align*}
where $L = (N^2-1)^{\frac{p-2}{2p-2}}/(2N)$. Since $f$ is convex, we have $f(x_k) - f(x_{k+1}) \ge \langle \nabla f(x_{k+1}), x_k-x_{k+1} \rangle$. Furthermore, since $f$ is $\sigma$-uniformly convex of order $p$, from~\cite[Lemma~3]{Nesterov08} we also have
\begin{align}\label{Eq:expGp3}
\|\nabla f(x_{k+1})\|_*^{\frac{p}{p-1}} \,\ge\, \frac{p}{p-1} \sigma^{\frac{1}{p-1}} (f(x_{k+1})-f(x^*)) \,\ge\, \sigma^{\frac{1}{p-1}} (f(x_{k+1})-f(x^*)).
\end{align}
Combining these inequalities and recalling the definition $\kappa = \epsilon \sigma$ gives us
\begin{align*}
f(x_k) - f(x_{k+1}) \,\ge\, L \kappa^{\frac{1}{p-1}} (f(x_{k+1})-f(x^*)),
\end{align*}
or equivalently,
\begin{align}\label{Eq:expGp2}
f(x_{k+1}) - f(x^*) \,\le\, \frac{f(x_k)-f(x^*)}{1+ L\kappa^{\frac{1}{p-1}}}
\,\le\, \frac{f(x_1)-f(x^*)}{\big(1+L \kappa^{\frac{1}{p-1}}\big)^k}.
\end{align}
Note that by the smoothness of $f$, as in~\eqref{Eq:AlgGradResidual2}, we can write $f(x_1) \le \min_x \{f(x) + \frac{N+1}{\epsilon p} \|x-x_0\|^p\} \le f(x^*) + \frac{N+1}{\epsilon p} \|x_0-x^*\|^p$.
Furthermore, since we assume the inverse condition number $\kappa = \epsilon \sigma$ is small, we can write $1+L \kappa^{\frac{1}{p-1}} \approx \exp(L \kappa^{\frac{1}{p-1}})$.
Therefore,~\eqref{Eq:expGp2} yields the desired convergence rate~\eqref{Eq:GradRateUnif}.
\end{proof}

Notice that the result of Theorem~\ref{Thm:GradRateUnif} matches the desired convergence rate $O(\frac{1}{\epsilon} e^{-ck})$ discussed in Appendix~\ref{App:Exp}, with $c = L \kappa^{\frac{1}{p-1}}$.

\paragraph{Exponential convergence rate of rescaled gradient flow.}
As a side remark, we note that the rescaled gradient flow also has an exponential convergence rate when the objective function $f$ is uniformly convex. However, notice that the following continuous-time convergence rate only depends on the uniform convexity constant of $f$, whereas the discrete-time convergence rate above also depends on the Lipschitz constant for the higher-order smoothness of $f$.

\begin{theorem}\label{Thm:RescUnif}
If $f$ is $\sigma$-uniformly convex of order $p$, then the rescaled gradient flow~\eqref{Eq:RescGradFlow} has convergence rate
\begin{align}\label{Eq:RescUnif}
f(X_t) - f(x^*) \le (f(X_0) - f(x^*)) \exp\left(-\sigma^{\frac{1}{p-1}} t\right).
\end{align}
\end{theorem}
\begin{proof}
As we saw in~\eqref{Eq:expGp3}, the uniform convexity of $f$ implies the inequality
\begin{align*}
\|\nabla f(X_t)\|_*^{\frac{p}{p-1}} \ge \sigma^{\frac{1}{p-1}} (f(X_t) - f(x^*)).
\end{align*}
Using this inequality and plugging in the rescaled gradient flow equation~\eqref{Eq:RescGradFlow}, we have
\begin{align*}
\frac{d}{dt} (f(X_t) - f(x^*)) \,&=\, \langle \nabla f(X_t), \dot X_t \rangle
  =\, -\|\nabla f(X_t)\|_*^{\frac{p}{p-1}} 
  \le\, -\sigma^{\frac{1}{p-1}} (f(X_t) - f(x^*)).
\end{align*}
Dividing both sides by $f(X_t) - f(x^*)$ and integrating, we get the desired convergence rate~\eqref{Eq:RescUnif}.
\end{proof}

\subsubsection{Exponential convergence rate of accelerated method with restart scheme}
\label{App:ExpConvergenceAcc}

We now show that a variant of the accelerated gradient method~\eqref{Eq:AlgComp} with a restart scheme also attains an exponential convergence rate, with a better dependence on the condition number $\kappa$ than the higher-order gradient method as in Appendix~\ref{App:ExpConvergenceGrad}. 

Specifically, we consider the following variant of the accelerated gradient method~\eqref{Eq:AlgComp}, 
\begin{subequations}\label{Eq:AlgComp2}
\begin{align}
x_{k+1} &= \frac{p}{k+p} z_k + \frac{k}{k+p} y_k   \label{Eq:AlgComp2x} \\
y_k &= \arg\min_y \left\{ f_{p-1}(y;x_k) + \frac{2}{\epsilon p} \|y-x_k\|^p \right\}  \label{Eq:AlgComp2y} \\
z_k &= \arg\min_z \left\{\frac{p}{(4p)^p} \sum_{i=0}^k i^{(p-1)} \langle \nabla f(y_i), z \rangle + \frac{2^{p-2}}{\epsilon p} \|z-x_0\|^p \right\}.  \label{Eq:AlgComp2z}
\end{align}
\end{subequations}
In~\eqref{Eq:AlgComp2}, for simplicity we have explicitly set the constant $N$ in~\eqref{Eq:AlgCompy} to be $N=2$, and set $C$ in~\eqref{Eq:AlgCompz} to be $C = 1/(4p)^p$, which satisfies the condition $C \le (N^2-1)^{\frac{p-2}{2}}/((2N)^{p-1} p^p)$. Furthermore, for the $z$-update~\eqref{Eq:AlgComp2z} we have used the equivalent version~\eqref{Eq:EstFunc0} where we unroll the recursion, and we have also replaced the Bregman divergence in the $z$-update~\eqref{Eq:AlgCompz} by the rescaled $p$-th power $d_p(z) = \frac{2^{p-2}}{p} \|z-x_0\|^p$, which is $1$-uniformly convex of order $p$. 
The proof of Theorem~\ref{Thm:AlgRate} still holds in this case, so we have the guarantee
\begin{align}\label{Eq:AlgRateMod}
f(y_k) - f(x^*) \,\le\, \frac{(4p)^p \cdot 2^{p-2} \|x_0-x^*\|^p}{\epsilon  p k^{(p)}}
\,\le\, \frac{2^{3p-2} p^{p-1} \|x_0-x^*\|^p}{\epsilon k^p}.
\end{align}

Then we define the following restart scheme, which proceeds by running the accelerated method~\eqref{Eq:AlgComp2} for some number of iterations at each step,
\begin{align}\label{Eq:expRest}
\hat x_{k} \,=\, \big(\text{the output $y_m$ of running~\eqref{Eq:AlgComp2} for $m$ iterations with input $x_0 = \hat x_{k-m}$}\big).
\end{align}
Our main result is the following.

\begin{theorem}\label{Thm:AccRateUnif}
Suppose $f$ is $\frac{(p-1)!}{\epsilon}$-smooth of order $p-1$, and $\sigma$-uniformly convex of order $p$. 
Let $\hat x_k$ be the output of running the restart scheme~\eqref{Eq:expRest} for $k/m$ times with $m = 8p/\kappa^{\frac{1}{p}}$, where $\kappa = \epsilon \sigma$ is the inverse condition number, and let $\hat y_k = G_{p,\epsilon,2}(\hat x_k)$ be the output of running one step of the gradient update~\eqref{Eq:Update} with input $\hat x_k$. Then we have the convergence rate
\begin{align}\label{Eq:expRestRate}
f(\hat y_k) - f(x^*) \,\le\, \frac{3\|\hat x_0-x^*\|^p}{\epsilon p \, e^{k/m}} \,=\, O\left(\frac{1}{\epsilon}  \, \exp\left(-\frac{\kappa^{\frac{1}{p}} k}{8 p}\right)\right).
\end{align}
\end{theorem}
\begin{proof}
Since $f$ is $\sigma$-uniformly convex of order $p$, and by the bound~\eqref{Eq:AlgRateMod}, we have
\begin{align}\label{Eq:expRest2}
\frac{\sigma}{p} \|\hat x_k-x^*\|^p
\,\le\, f(\hat x_k)-f(x^*)
\,\le\, \frac{2^{3p-2} p^{p-1} \|\hat x_{k-m}-x^*\|^p}{\epsilon m^p}
\,\le\, \frac{\sigma}{pe} \|\hat x_{k-m}-x^*\|^p,
\end{align}
where the last inequality follows from our choice of $m$. 
Thus, an execution of~\eqref{Eq:expRest} with $m$ iterations of the accelerated method reduces the distance to optimum by a factor of at least $1/e$. 
Iterating~\eqref{Eq:expRest2}, we obtain $\|\hat x_k-x^*\|^p \le e^{-k/m}\|\hat x_0-x^*\|^p$.
To convert this into a bound on the function value, we use the smoothness of $f$. As noted in~\eqref{Eq:AlgGradResidual2}, since $\hat y_k$ is the output of one step of the gradient update~\eqref{Eq:Update} with input $\hat x_k$, we have $f(\hat y_k) - f(x^*) \le \frac{3}{\epsilon p} \|\hat x_k-x^*\|^p$. 
This gives the desired bound~\eqref{Eq:expRestRate}.
\end{proof}

The result of Theorem~\ref{Thm:AccRateUnif} matches the desired convergence rate $O(\frac{1}{\epsilon} e^{-ck})$ as discussed in Appendix~\ref{App:Exp} with $c = \frac{1}{8p} \kappa^{\frac{1}{p}}$. Note that this convergence rate has a better dependence on the inverse condition number $\kappa = \epsilon \sigma$ than the higher-order gradient algorithm as in Theorem~\ref{Thm:GradRateUnif}, because $\kappa^{\frac{1}{p}} > \kappa^{\frac{1}{p-1}}$ for small $\kappa$. This generalizes the conclusion of~\cite[Section~5]{Nesterov08} for the case $p=3$.
However, as noted previously, the link to continuous time is not as clear as that of the polynomial family.

\subsection{Hessian vs.\ Bregman Lagrangian}
\label{App:HessLag}

In a Hessian manifold, the metric is generated by the Hessian $\nabla^2 h$ of the distance-generating function $h$.
So for example, the gradient flow equation in the Euclidean case, $\dot X_t = -\nabla f(X_t)$, which can be written as
\begin{align*}
\dot X_t = \arg\min_V \left\{ \langle \nabla f(X_t), V \rangle + \frac{1}{2} \|V\|^2 \right\},
\end{align*}
in general becomes the natural gradient flow $\dot X_t = -[\nabla^2 h(X_t)]^{-1} \nabla f(X_t)$, or equivalently,
\begin{align*}
\dot X_t = \arg\min_V \left\{ \langle \nabla f(X_t), V \rangle + \frac{1}{2} \|V\|^2_{\nabla^2 h(X_t)} \right\},
\end{align*}
which is obtained by replacing the Euclidean squared norm $\|v\|^2 = \langle v,v \rangle$ by the Hessian metric
\begin{align*}
\|v\|^2_{\nabla^2 h(x)} := \langle \nabla^2 h(x) v, v \rangle.
\end{align*}

At the Lagrangian level, recall that a starting point of our work is the differential equation $\ddot X_t + \frac{3}{t} \dot X_t + \nabla f(X_t) = 0$ for accelerated gradient descent~\cite{SuBoydCandes14}, which we observe is the Euler-Lagrange equation for the damped Lagrangian
\begin{align}\label{Eq:LagSBC}
\L(X,V, t) = t^3\left(\frac{1}{2} \|V\|^2 - f(X)\right).
\end{align}
How should we generalize this Lagrangian to the non-Euclidean case? 
From our discussion on natural gradient flow, a natural guess is to replace the Euclidean metric in~\eqref{Eq:LagSBC} by the Hessian metric. 
Thus, we are led to consider the following family of \emph{Hessian Lagrangians}:
\begin{align}\label{Eq:HessLag}
\L_{\text{Hess}}(X,V,t) \,=\, e^{\gamma_t} \left(\frac{1}{2} \|V\|^2_{\nabla^2 h(X)} - e^{\beta_t} f(X) \right)
\end{align}
where we have also introduced arbitrary weighting functions $\beta_t, \gamma_t \in \R$ (\eqref{Eq:LagSBC} is the Euclidean case with $\beta_t = 0, \gamma_t = 3 \log t$).
However, the Hessian Lagrangian~\eqref{Eq:HessLag} turns out to be unsuitable for our optimization purposes.
This is because the Euler-Lagrange equation for the Hessian Lagrangian~\eqref{Eq:HessLag},
\begin{align}\label{Eq:ELHess}
\frac{1}{2} \nabla^3 h(X_t) \, \dot X_t \, \dot X_t \,+\, \nabla^2 h(X_t) \left(\ddot X_t + \dot \gamma_t \, \dot X_t \right) + e^{\beta_t}\nabla f(X_t) = 0,
\end{align}
involves the third-order derivative $\nabla^3 h$ (which comes from being the derivative of the metric tensor $\nabla^2 h$). This makes the analysis difficult, preventing us from obtaining a convergence rate for~\eqref{Eq:ELHess}. Furthermore, the presence of $\nabla^3 h$ in the equation makes it difficult to implement as an efficient discrete-time algorithm.

On the other hand, our work shows that the ``correct'' way to generalize~\eqref{Eq:LagSBC} to the non-Euclidean case is to use the Bregman divergence, rather than Hessian metric. This results in the general Bregman Lagrangian family~\eqref{Eq:BregLag}, which requires an additional parameter $\alpha_t$ controlling the amount of interaction between the position $X$ and velocity $V$. When the parameters are coupled in an ideal scaling, the Bregman Lagrangian produces dynamics that converge at a provable rate. This is achieved via the design of a corresponding Lyapunov function (the energy functional $\E_t$~\eqref{Eq:E}), whose form is intimately tied to the use of the Bregman divergence in the Lagrangian. Furthermore, for the polynomial family, we can discretize the resulting dynamics as a discrete-time algorithm~\eqref{Eq:AlgComp} that does not require the Hessian $\nabla^2 h$, but only the gradient $\nabla h$.

It is interesting to consider whether the Hessian Lagrangian~\eqref{Eq:HessLag} has useful properties, and how it relates to the Bregman Lagrangian.
For a small displacement $\varepsilon > 0$ we know that Bregman divergence approximates the Hessian metric, i.e., $D(x + \varepsilon v, x) \approx \frac{\varepsilon^2}{2} \|v\|_{\nabla^2 h(x)}^2$. Setting $\varepsilon = e^{-\alpha_t}$, this suggests that the Bregman Lagrangian~\eqref{Eq:BregLag} is approximating the Hessian Lagrangian $\L_{\text{Hess}}(X,V, t) =  e^{\gamma_t-\alpha_t} \left(\frac{1}{2} \|V\|_{\nabla^2 h(X)}^2 - e^{2\alpha_t+\beta_t} f(X) \right)$.
However, this argument assumes $\epsilon$ is small, whereas in our particular case of interest (the polynomial subfamily in Section~\ref{Sec:Poly}) the value of $\epsilon = e^{-\alpha_t} = \frac{t}{p}$ is growing over time.

\subsection{Gradient vs.\ Lagrangian flows}
\label{App:GradLag}
  
  In the Euclidean case, we can think of gradient flow as describing the behavior of a damped Lagrangian system ``in an asymptotic regime in which dissipative effects play such an important role, that the effects of forcing and dissipation compensate each other''~\cite[p.~646]{Villani}.
  That is, the gradient flow equation $\dot X_t = -\nabla f(X_t)$ can be seen as the strong-friction limit $\lambda \to \infty$ of the equation $\ddot X_t + \lambda \dot X_t + \lambda \nabla f(X_t) = 0$. 
  
  This is perhaps more apparent if we define $m = 1/\lambda$ to be the ``mass'' of the fictitious particle, so the equation of motion becomes
  \begin{align}\label{Eq:m0}
  m\ddot X_t + \dot X_t + \nabla f(X_t) = 0,
  \end{align}
  which is the Euler-Lagrange equation of the damped Lagrangian
  \begin{align}\label{Eq:m0Lag}
  \L(X,V,t) = e^{t/m}\left(\frac{m}{2} \|V\|^2 - f(X)\right),
  \end{align}
  where the damping factor $e^{t/m}$ also scales with $m$. 
  In the massless limit $m \to 0$, we indeed recover gradient flow from~\eqref{Eq:m0}. 
  In the following, we show that this result also holds more generally, both for natural gradient flow (as the massless limit of a Bregman Lagrangian flow) and for the rescaled gradient flow (as the massless limit of a Lagrangian flow which uses the $p$-th power of the norm). 

  However, notice that in all these cases, the momentum variable $P = \frac{\partial \L}{\partial V}$ becomes infinite as $m \to 0$. For instance, $P = me^{t/m} V$ for~\eqref{Eq:m0Lag}, and $me^{t/m} \to \infty$. This means as $m \to 0$, the particle also becomes more massive and has more inertia. Thus, gradient flow is the limiting case where the infinitely massive particle simply rolls downhill and stops at the minimum $x^*$ as soon as the force $-\nabla f$ vanishes, without oscillation (which is damped by the infinitely strong friction).
  In this view, moving from a first-order gradient algorithm to a second-order Lagrangian (accelerated) algorithm does {\em not} amount to preventing oscillation; rather, it is the opposite, by unwinding the curve to finite momentum where it can travel faster, albeit with some oscillation.

\paragraph{Natural gradient flow as massless limit.}
Consider the following Lagrangian
\begin{align}\label{Eq:mBregLag}
\L(X,V,t) = \frac{e^{t/m}}{m} \left(D_h(X + mV, X) - mf(X) \right),
\end{align}
which is the Bregman Lagrangian~\eqref{Eq:BregLag} with parameters $\alpha_t = -\log m$, $\beta_t = \log m$, and $\gamma_t = t/m$ (which satisfy the ideal scaling~\eqref{Eq:IdeSca}). Note that~\eqref{Eq:mBregLag} recovers~\eqref{Eq:m0Lag} in the Euclidean case. 
The Euler-Lagrange equation~\eqref{Eq:ELBreg2} for the Lagrangian~\eqref{Eq:mBregLag} is given by
\begin{align*}
\ddot X_t + \frac{1}{m} \dot X_t + \frac{1}{m} \left[\nabla^2 h(X_t + m \dot X_t)\right]^{-1} \nabla f(X_t) = 0.
\end{align*}
Multiplying the equation by $m$ and letting $m \to 0$, we recover
\begin{align*}
\dot X_t + \left[\nabla^2 h(X_t)\right]^{-1} \nabla f(X_t) = 0,
\end{align*}
which is the natural gradient flow equation.
In this case the momentum variable is $P = \frac{\partial \L}{\partial V} = e^{t/m}(\nabla h(X + m V) - \nabla h(X)) \approx me^{t/m} \nabla^2 h(X) V$, so we still have $P \to \infty$ as $m \to 0$.

\paragraph{Rescaled gradient flow as massless limit.}
Consider the following Lagrangian
\begin{align}\label{Eq:pLag}
\L(X,V,t) = e^{t/m} \left(\frac{m}{p} \|V\|^p - f(X)\right),
\end{align}
where we use the $p$-th power of the norm to measure the kinetic energy. Note that~\eqref{Eq:pLag} recovers~\eqref{Eq:m0Lag} in the case $p=2$. The Euler-Lagrange equation is
\begin{align*}
\|\dot X_t\|^{p-2} \big(m\ddot X_t + \dot X_t \big) + m(p-2) \|\dot X_t\|^{p-4} \langle \ddot X_t,\dot X_t \rangle \, \dot X_t + \nabla f(X_t) = 0.
\end{align*}
So as $m \to 0$, this equation recovers
\begin{align}
\|\dot X_t\|^{p-2} \dot X_t + \nabla f(X_t) = 0
\end{align}
which is equivalent to the rescaled gradient flow~\eqref{Eq:RescGradFlow}.
In this case the momentum variable is $P = \frac{\partial \L}{\partial V} = me^{t/m}\|V\|^{p-2} V$, which still goes to infinity as $m \to 0$.

\subsection{Bregman Hamiltonian}
\label{App:BregHam}

In this section we define and compute the Bregman Hamiltonian corresponding to the Bregman Lagrangian.
In general, given a Lagrangian $\L(X,V,t)$, its Hamiltonian is defined by
\begin{align}\label{Eq:HamDef}
\H(X,P,t) = \langle P,V \rangle - \L(X, V, t)
\end{align}
where $P = \frac{\partial \L}{\partial V}$ is the momentum variable conjugate to position.

For the Bregman Lagrangian~\eqref{Eq:BregLag}, the momentum variable is given by
\begin{align}\label{Eq:BregMom}
P =  \frac{\partial \L}{\partial V} = e^{\gamma_t} \left( \nabla h(X + e^{-\alpha_t} V) - \nabla h(X)  \right).
\end{align}
We can invert this equation to solve for the velocity $V$,
\begin{align}\label{Eq:BregxDot}
V \,=\, e^{\alpha_t} \left( \nabla h^* (\nabla h(X) + e^{-\gamma_t} P) - X \right),
\end{align}
where $h^*$ is the conjugate function to $h$ (recall the definition in~\eqref{Eq:hDual}), and we have used the property that $\nabla h^* = [\nabla h]^{-1}$. So for the first term in the definition~\eqref{Eq:HamDef} we have
\begin{align*}
\langle P, V \rangle = e^{\alpha_t} \big\langle P, \; \nabla h^* (\nabla h(X) + e^{-\gamma_t} P) - X  \big\rangle.
\end{align*}

Next, we write the Bregman Lagrangian $\L(X,V,t)$ in terms of $(X,P,t)$. We can directly substitute~\eqref{Eq:BregxDot} to the definition~\eqref{Eq:BregLag} and calculate the result. Alternatively, we can use the property that the Bregman divergences of $h$ and $h^*$ satisfy $D_h(y,x) = D_{h^*}(\nabla h(x), \nabla h(y))$.
Therefore, we can write the Bregman Lagrangian~\eqref{Eq:BregLag} as
\begin{align*}
\L(X,V, t)
    &= e^{\alpha_t+\gamma_t} \left(D_{h^*}\left(\nabla h(X), \, \nabla h(X + e^{-\alpha_t} V)\right) - e^{\beta_t} f(X) \right)   \\
    &= e^{\alpha_t+\gamma_t} \left(D_{h^*}\left(\nabla h(X), \, \nabla h(X) + e^{-\gamma_t} P\right) - e^{\beta_t} f(X) \right)  \\
    &= \, e^{\alpha_t+\gamma_t} \left(h^*(\nabla h(X)) - h^*(\nabla h(X) + e^{-\gamma_t} P) + e^{-\gamma_t} \langle \nabla h^*(\nabla h(X) + e^{-\gamma_t} P), P \rangle - e^{\beta_t} f(X) \right),
\end{align*}
where in the second step we have used the relation $\nabla h(X+e^{-\alpha_t} V) = \nabla h(X) + e^{-\gamma_t} P$ from~\eqref{Eq:BregMom}, and in the last step we have expanded the Bregman divergence.

Substituting these calculations into~\eqref{Eq:HamDef} and simplifying, we get the Hamiltonian
\begin{align*}
\H(X,P,t) \,=\, e^{\alpha_t+\gamma_t} \left(h^*(\nabla h(X) + e^{-\gamma_t} P) - h^*(\nabla h(X)) - \langle X, \, e^{-\gamma_t} P \rangle + e^{\beta_t} f(X) \right).
\end{align*}
Since $X = \nabla h^*(\nabla h(X))$, we can also write this result in terms of the Bregman divergence of $h^*$,
\begin{align}\label{Eq:BregHam}
\H(X,P,t) \,=\, e^{\alpha_t+\gamma_t} \left(D_{h^*}(\nabla h(X) + e^{-\gamma_t} P, \, \nabla h(X)) + e^{\beta_t} f(X) \right).
\end{align}
We call the Hamiltonian~\eqref{Eq:BregHam} the \emph{Bregman Hamiltonian}. Notice that whereas the Bregman Lagrangian takes the form of the difference between the kinetic and potential energy, the Bregman Hamiltonian takes the form of the sum of the kinetic and potential energy. (However, note that the kinetic energy is slightly different: it is $D_{h^*}(\nabla h(X) + e^{-\gamma_t} P, \, \nabla h(X)) = D_h(X, X+e^{-\alpha_t} V)$ in the Hamiltonian~\eqref{Eq:BregHam}, while it is $D_h(X+e^{-\alpha_t} V, X)$ in the Lagrangian~\eqref{Eq:BregLag}.)

\paragraph{Hamiltonian equations of motion.}
The second-order Euler-Lagrange equation of a Lagrangian can be equivalently written as a pair of first-order equations
\begin{equation}\label{Eq:HamFlow}
\dot X_t \,=\, \frac{\partial \H}{\partial P} (X_t,P_t,t),
\quad\quad
\dot P_t \,=\, -\frac{\partial \H}{\partial X} (X_t,P_t,t).
\end{equation}

For the Bregman Hamiltonian~\eqref{Eq:BregHam}, the equations of motion are given by
\begin{subequations}\label{Eq:HamEq}
\begin{align}
\dot X_t \,&=\, e^{\alpha_t} \left( \nabla h^*(\nabla h(X_t) + e^{-\gamma_t} P_t) - X_t \right)   \label{Eq:HamEqx} \\
\dot P_t \,&=\, -e^{\alpha_t+\gamma_t} \nabla^2 h(X_t) \left( \nabla h^*(\nabla h(X_t) + e^{-\gamma_t} P_t) - X_t \right) + e^{\alpha_t} P_t - e^{\alpha_t+\beta_t+\gamma_t} \nabla f(X_t)  \label{Eq:HamEqp}.
\end{align}
\end{subequations}
Notice that the first equation~\eqref{Eq:HamEqx} recovers the definition of momentum~\eqref{Eq:BregMom}. Furthermore, when $\dot \gamma_t = e^{\alpha_t}$, by substituting~\eqref{Eq:HamEqx} to~\eqref{Eq:HamEqp} we can write~\eqref{Eq:HamEq} as
\begin{align*}
\frac{d}{dt} \left\{\nabla h(X_t) + e^{-\gamma_t} P_t \right\}
   \,=\, \nabla^2 h(X_t) \, \dot X_t - \dot \gamma_t e^{-\gamma_t} P_t + e^{-\gamma_t} \dot P_t
    \,=\, - e^{\alpha_t+\beta_t} \nabla f(X_t).
\end{align*}
Since $\nabla h(X_t) + e^{-\gamma_t} P_t = \nabla h(X_t + e^{-\alpha_t} \dot X_t)$ by~\eqref{Eq:HamEqx}, this indeed recovers the Euler-Lagrange equation~\eqref{Eq:ELBreg3}.

A Lyapunov function for the Hamiltonian equations of motion~\eqref{Eq:HamEq} is the following, which is simply the energy functional~\eqref{Eq:E} written in terms of $(X_t,P_t,t)$,
\begin{align*}
\E_t \,=\, D_{h^*}\left(\nabla h(X_t) + e^{-\gamma_t} P_t, \, \nabla h(x^*)\right) + e^{\beta_t} (f(X_t)-f(x^*)).
\end{align*}

The Hamiltonian formulation of the dynamics has appealing properties that seem worthy of further exploration. For example, 
Hamiltonian flow preserves volume in phase space (Liouville's theorem); this property has been used in the context of sampling to develop the technique of Hamiltonian Markov chain Monte-Carlo, and may also be useful to help us design better algorithms for optimization.
Furthermore, the Hamilton-Jacobi-Bellman equation (which is a reformulation of the Hamiltonian dynamics) is a central object of study in the field of optimal control theory, and it would be interesting to study the Bregman Hamiltonian framework from that perspective.

\subsection{Gauge invariance}
\label{App:Gauge}

The Euler-Lagrange equation of a Lagrangian is gauge-invariant, which means it does not change when we transform the Lagrangian by adding a total time derivative,
\begin{align}\label{Eq:LagT}
\L'(X_t,\dot X_t,t) = \L(X_t,\dot X_t,t) + \frac{d}{dt} G(X_t,t)
\end{align}
for any smooth function $G$. We can show this by directly checking that the Euler-Lagrange equation of $\L'$ is the same as that of $\L$. Alternatively, this follows from the formulation of the principle of least action, where we fix two points $(x_0,t_0)$ and $(x_1,t_1)$, and ask for a curve $X$ joining the two endpoints ($X_{t_0} = x_0$ and $X_{t_1} = x_1$) that minimizes the action $J(X) = \int_{t_0}^{t_1} \L(X_t, \dot X_t, t) dt$.
Thus, when the Lagrangian transforms as~\eqref{Eq:LagT}, the action only changes to $J'(X) = J(X) + \int_{t_0}^{t_1} \frac{d}{dt} G(X_t, t) dt = J(X) + G(x_1,t_1) - G(x_0,t_0)$.
Since $(x_0,t_0)$ and $(x_1,t_1)$ are fixed, this means the new action only differs from the old action by a constant; this implies that the optimal least action curve---namely, the Euler-Lagrange equation---does not change.

In our case, under the ideal scaling condition $\dot \gamma_t = e^{\alpha_t}$~\eqref{Eq:IdeScaGam}, this property implies that the Bregman Lagrangian~\eqref{Eq:BregLag} is equivalent to the following Lagrangian
\begin{align}\label{Eq:BregLagRed}
\L'(X,V,t) = e^{\gamma_t + \alpha_t} \left(h(X + e^{-\alpha_t} V) - e^{\beta_t} f(X) \right),
\end{align}
where we have replaced the Bregman divergence $D_h(X+e^{-\alpha_t} V, X)$ by its first term $h(X+e^{-\alpha_t} V)$.
Indeed, we can check that the difference between the Bregman Lagrangian~\eqref{Eq:BregLag} and the reduced form~\eqref{Eq:BregLagRed} is a total time derivative,
\begin{align*}
\L'(X_t,\dot X_t,t) - \L(X_t,\dot X_t,t) = e^{\gamma_t + \alpha_t} \left( h(X_t) + \langle \nabla h(X_t), e^{-\alpha_t} \dot X_t \rangle \right)
\,=\, \frac{d}{dt} \left\{ e^{\gamma_t} h(X_t) \right\}
\end{align*}
where the last step follows from the ideal scaling $e^{\alpha_t} = \dot \gamma_t$.

The reduced Lagrangian~\eqref{Eq:BregLagRed} is slightly simpler than the Bregman Lagrangian~\eqref{Eq:BregLag}, and in a sense it makes the roles of $h$ and $f$ more symmetric.
It also suggests that the role of $h$ is not so much as measuring the distance via the Hessian metric or Bregman divergence, but rather, as evaluating the extrapolated future point $X_t + e^{-\alpha_t} \dot X_t$.

\subsection{Natural motion}
\label{App:NatMot}

A natural motion is the motion of a particle when it experiences no force. In the physical world, the natural motion of a particle is a straight-line motion with constant velocity. But for the Bregman Lagrangian, which describes a dissipative system, the natural motion always converges.

Specifically, the Bregman Lagrangian~\eqref{Eq:BregLag} in the case of zero (or constant) potential function $f \equiv 0$ is $\L(X,V,t) = e^{\alpha_t + \gamma_t} D_h(X + e^{-\alpha_t} V, X)$. 
Assuming the ideal scaling $\dot \gamma_t = e^{\alpha_t}$~\eqref{Eq:IdeScaGam}, its Euler-Lagrange equation is given by~\eqref{Eq:ELBreg3}, which in this case is
\begin{align}\label{Eq:ELNat}
\frac{d}{dt} \nabla h(X_t + e^{-\alpha_t} \dot X_t) = 0.
\end{align}
This means $\nabla h(X_t + e^{-\alpha_t} \dot X_t)$ is a constant, say $\nabla h(X_t + e^{-\alpha_t} \dot X_t) = \nabla h(b)$ for some $b \in \X$. 
Applying $\nabla h^* = [\nabla h]^{-1}$ to both sides gives us $X_t + e^{-\alpha_t} \dot X_t = b$. 
Since $e^{\alpha_t} = \dot \gamma_t$, we can write this as
\begin{align*}
\frac{d}{dt} \left\{ e^{\gamma_t} (X_t-b) \right\} 
\,=\, e^{\alpha_t + \gamma_t} (X_t-b) + e^{\gamma_t} \dot X_t
\,=\, 0.
\end{align*}
This means $e^{\gamma_t}(X_t - b)$ is a constant, say $e^{\gamma_t}(X_t - b) = a$ for some $a \in \X$. Thus, we conclude that the natural motion of the Bregman Lagrangian is
\begin{align}\label{Eq:ELNat2}
X_t = ae^{-\gamma_t} + b.
\end{align}
Notice that the natural motion is independent of $h$, although the Lagrangian still depends on $h$.
Furthermore, in contrast with the straight-line motion, the natural motion~\eqref{Eq:ELNat2} always converges; in particular, if we assume $e^{\gamma_t} \to \infty$ as $t \to \infty$, then $X_t \to b$.

The natural motion~\eqref{Eq:ELNat2} has simple explicit invariance and symmetry properties. Indeed,~\eqref{Eq:ELNat} states that $\nabla h(X_t + e^{-\alpha_t} \dot X_t)$ is a conserved quantity, which is always equal to $\nabla h(b)$. By Noether's theorem, any conservation law corresponds to a symmetry of the Lagrangian. In our case, the corresponding symmetry is the transformation
\begin{align}\label{Eq:SymTrans}
X_t' = X_t + e^{-\gamma_t} u,
\end{align}
for any $u \in \X$. Under this transformation, $\dot X_t$ changes to $\dot X_t' = \dot X_t - \dot \gamma_t e^{-\gamma_t} u$. Since $\dot \gamma_t = e^{\alpha_t}$, this implies $X_t' + e^{-\alpha_t} \dot X_t' = X_t + e^{-\alpha_t} \dot X_t$. This means the reduced Lagrangian $\L(X_t,\dot X_t,t) = e^{\gamma_t + \alpha_t} h(X_t + e^{-\alpha_t} \dot X_t)$ is invariant under the transformation~\eqref{Eq:SymTrans}. Therefore, the Bregman Lagrangian (which is gauge-equivalent to the reduced Lagrangian) is also invariant. So indeed~\eqref{Eq:SymTrans} is a symmetry of the Bregman Lagrangian when $f = 0$.

\subsection{The Euclidean case}
\label{App:Euc}

In the Euclidean case many of our results and equations simplify, as we summarize in this section.
When $h$ is the squared Euclidean norm, $h(x) = \frac{1}{2} \|x\|^2$, the Bregman divergence is also the squared norm and it coincides with the Hessian metric, $D_h(y,x) = \frac{1}{2} \|y-x\|^2 = \frac{1}{2} \|y-x\|^2_{\nabla^2 h(x)}$. Furthermore, $h^* = h$ and both $\nabla h, \nabla h^*$ are the identity function.

In the Euclidean case, the Bregman Lagrangian~\eqref{Eq:BregLag} becomes
\begin{align*}
\L(X,V,t) = e^{\gamma_t - \alpha_t} \left(\frac{1}{2} \|V\|^2 - e^{2\alpha_t + \beta_t} f(X) \right).
\end{align*}
For general $\alpha_t, \beta_t, \gamma_t$, the Euler-Lagrange equation~\eqref{Eq:ELBreg} is given by
\begin{align*}
\ddot X_t + (\dot \gamma_t - \dot \alpha_t) \dot X_t + e^{2\alpha_t + \beta_t} \nabla f(X_t) = 0.
\end{align*}
When the ideal scaling $\dot \gamma_t = e^{\alpha_t}$~\eqref{Eq:IdeScaGam} holds, this equation becomes
\begin{align}\label{Eq:ELEuc}
\ddot X_t + (e^{\alpha_t} - \dot \alpha_t) \dot X_t + e^{2\alpha_t + \beta_t} \nabla f(X_t) = 0,
\end{align}
which we can equivalently write as $\frac{d}{dt} (X_t + e^{-\alpha_t} \dot X_t ) = -e^{\alpha_t + \beta_t} \nabla f(X_t)$.
The energy functional~\eqref{Eq:E} for proving the rate of convergence becomes
\begin{align*}
\E_t = \frac{1}{2} \|X_t + e^{-\alpha_t} \dot X_t - x^*\|^2 + e^{\beta_t}(f(X_t) - f(x^*)).
\end{align*}
The Bregman Hamiltonian~\eqref{Eq:BregHam} becomes
\begin{align*}
\H(X,P,t) = e^{\alpha_t - \gamma_t} \left( \frac{1}{2} \|P\|^2 + e^{2\gamma_t + \beta_t} f(X) \right),
\end{align*}
where the momentum variable~\eqref{Eq:BregMom} is given by $P = e^{\gamma_t - \alpha_t} V$.
The Hamiltonian equations of motion~\eqref{Eq:HamEq} simplify to
\begin{align}
\dot X_t \,&=\, e^{\alpha_t - \gamma_t} P_t \\
\dot P_t \,&=\, - e^{\alpha_t+\beta_t+\gamma_t} \nabla f(X_t).
\end{align}

In particular, for the polynomial case with the parameters~\eqref{Eq:ABCPoly}, 
the Euler-Lagrange equation~\eqref{Eq:ELEuc} is given by
\begin{align}\label{Eq:ELEucPoly}
\ddot X_t + \frac{p+1}{t} \dot X_t + Cp^2 t^{p-2} \nabla f(X_t) = 0,
\end{align}
with an $O(1/t^p)$ rate of convergence.
For $p=2$, this recovers the differential equation $\ddot X_t + \frac{3}{t} + \nabla f(X_t) = 0$ corresponding to Nesterov's accelerated gradient descent, as derived in~\cite{SuBoydCandes14}. 

Su et al.~\cite{SuBoydCandes14} observed that the generalized equation $\ddot X_t + \frac{r}{t} \dot X_t + \nabla f(X_t) = 0$ still has convergence rate $O(1/t^2)$ whenever $r \ge 3$, and they posed the question on the significance of the threshold $r = 3$.
Our results give the following perspective: The equation $\ddot X_t + \frac{r}{t} \dot X_t + \nabla f(X_t) = 0$ is the case of~\eqref{Eq:ELEuc} with parameters $\alpha_t = \log (r-1) - \log t$, $\gamma_t = (r-1) \log t$, and $\beta_t = 2\log t - 2\log (r-1)$. These parameters satisfy the ideal scaling condition~\eqref{Eq:IdeSca} when $r \ge 3$, so Theorem~\ref{Thm:Rate} guarantees a convergence rate of $O(e^{-\beta_t}) = O(1/t^2)$. 
However, for a fixed $r > 3$, the choice of $\beta_t = 2\log t - 2\log (r-1)$ is suboptimal, since from the ideal scaling condition $\dot \beta_t \le e^{\alpha_t}$ we know we can increase $\beta_t$ up to $(r-1) \log t$. This will introduce a factor of $t^{r-3}$ on the force term, as in~\eqref{Eq:ELEucPoly}, but it will also yield a faster convergence rate of $O(1/t^{r-1})$.

\end{document}